\documentclass[12pt]{amsart}

\usepackage[textsize=footnotesize,backgroundcolor=yellow!70,bordercolor=orange]{todonotes}
\definecolor{refkey}{rgb}{.35,.75,0}
\definecolor{labelkey}{rgb}{.15,.55,0}

\usepackage{filemod}
\usepackage{graphicx}
\graphicspath{{../graf/}{./}{./graf_new/}}
\usepackage{soul}
\usepackage{psfrag}
\usepackage{stmaryrd}
\usepackage{cancel}
\usepackage{multirow}
\usepackage{xcolor} 
\usepackage{amsmath} 
\usepackage{amssymb}
\usepackage[percent]{overpic}
\usepackage{amsfonts, dsfont} 
\usepackage{subfig}

\usepackage{amsthm, amsfonts}

\newtheorem{example}{Example}
\newtheorem{remark}{Remark}

\usepackage{nicefrac}

\usepackage{tikz}
\usetikzlibrary{shapes}

\newcommand{\R}{\mathbb{R}}

\def\maxd{\mathop{{\max}}\limits}

\newcommand{\ds}{\displaystyle}

\title[Coupling conditions for shallow water junctions]{Angle dependence  in coupling conditions for shallow water equations at canal junctions}

\author{M. Briani $^{1}$, G. Puppo $^{2}$, M. Ribot $^{3}$}

\begin{document}

\footnote{Istituto per le Applicazioni del Calcolo, Consiglio Nazionale delle Ricerche, Rome, Italy (m.briani@iac.cnr.it)}
\footnote{Sapienza Universit\`a, Rome, Italy (gabriella.puppo@uniroma1.it)}
\footnote{Universit\'e d'Orl\'eans, France (magali.ribot@univ-orleans.fr)}
 
\maketitle

\begin{abstract} 
In this paper we propose a numerical Riemann problem solver at the junction of one dimensional shallow-water canal networks. The junction conditions take into account the angles with which the channels intersect and include the possibility of canals with different sections. The solver is illustrated with several numerical tests which underline the importance of the angle dependence to obtain reliable solutions.  
\end{abstract}

\section{Introduction}
The shallow water equations are a model to describe free surface water flows.
They are a non-linear hyperbolic system of PDEs consisting of a mass and momentum balance. They are used to describe flows in artificial canals and water channels with applications for instance to environmental problems.
In water management problems, these equations are often used as a fundamental tool to describe the dynamics of networks of canals or of the branching of rivers. Networks occur in different type of configurations. 
The most straightforward treatment from a numerical point of view consists in considering the network as a two dimensional domain covered with an unstructured grid \cite{FullSWOFpaper, DelisKatsounis2D}.
However from a computational point of view it is much more efficient to consider the network as a set of one dimensional canals coupled through junctions.

The main difficulty mathematically is the definition of
 the coupling conditions at the junction between the adjoining channels.
We mention the following reviews for one-dimensional flows on networks \cite{BCGHP14,garavello2010review}. 
The coupling condition can be seen as a Riemann problem involving a constant state for each of the adjoining canals. Riemann problems at a junction are widely discussed in literature, see \cite{GP09,CHS2008,Kreiss2013,Marigo2010,LS2002}.
To close the problem one completes the Riemann problem with physical conservation properties across the junction, see \cite{Kesserwani2008,Jacovkis1991}.

From a numerical point of view, one has to couple one dimensional numerical solvers in the 1D channels with an approximate junction Riemann solver, see \cite{Toro2018,BorscheNetwork,BPQ2016,HS2013,roggensack2013kinetic}.

We consider a junction of three canals and we assume by convention to
have one incoming canal which ends at the junction and two outgoing canals which start at the junction. To solve the junction problem we need to find 
the three states (mass and discharge) facing each of the three one dimensional channels at the junction for a total of six unknowns.
One imposes mass conservation at the junction which yields one equation, 
then one formulates a left-half Riemann problem for the incoming
canal and a right-half Riemann problem for each of the two outgoing canals. Under subcritical flow assumptions, we obtain therefore three more conditions. Thus, two remaining equations have to be specified in order to define the junction model. In some works the
set of equations is completed by assuming the continuity of water levels \cite{CHS2008,Kreiss2013,Marigo2010,BP2018} or the continuity of energy \cite{Kesserwani2008,HS2013}. 
However, none of these works use a condition that takes into account the geometry and especially the angles formed by the channels in the fork. 
For an attempt to include an angle dependency in the solver, see \cite{CDG2019,HR, Ghostine2009}. 
Different approaches covert the junction with two-dimensional elements and project the computed 2D solution along the one dimensional channels \cite{Toro2018}.

These studies have several applications, such as optimization \cite{LS2002,GLS2004,G2005,GL09} or see \cite{bocchi2019} for a nice application to the modeling of a particular wave energy converter, the so-called oscillating water column.

In this article, we propose new coupling conditions at the junction that depend on the angles with which the channels intersect at the junction allowing also for channels with different sections, \cite{EscalanteSemplice1D, LeflochThanh1D}. 
Away from the junction we assume the solution to be 1D, while we describe the junction as a 2D region where coupling occurs between the branches. 
 We then consider the triangle formed by the intersection points of the walls of the three channels and, to this two dimensional domain, we apply conservation of mass and of the two components of momentum.  We obtain three non linear equations which include a dependence on the angles, for the six unknowns at the junction to be coupled with the three equations of the characteristic curves. Extending this study to a network of canals with several nodes is straightforward. This work therefor extends the results of \cite{CDG2019} by considering  branches with different sections. 
 
We prove that the conservation condition based on the continuity of the
energy across the junction occurs in a particular configuration of our setting. We also prove the existence of the solution of our junction Riemann problem in a few particular cases. 

Validation of numerical schemes obtained in this way is carried out comparing the numerical 1D solution with the junction, with a fully 2D solver, see 
\cite{Ghostine2009, Ghostine2013, HNS2019, Ghostine2012}.
We compare our numerical solver with a fully 2D solver for shallow water equations showing that the numerical approximation improves as the width of the 2D channels is reduced.

The paper is organized as follows. In Section \ref{sec:sw_Rp}, we concentrate on the solution of the Riemann problem for shallow-water equations. We then present the 
junction geometry in Section \ref{sec:Junction}, defining our coupling conditions in Section \ref{sec:Coupling_junction}. We discuss extensions for special configurations in Section \ref{sec:j_special}, this includes the case of a single channel with varying cross-section. In Section \ref{sec:NumScheme} we merge the relations at the junction with the numerical approximation of shallow water equations along the channels. Section \ref{sec:Sol_System} is devoted to a discussion of the existence of the numerical solution in a few cases.
We end in Section \ref{Sec:NumTests} with the numerical tests.

\section{The shallow water equations and its standard Riemann problem}\label{sec:sw_Rp}

Let us first recall the shallow water or Saint Venant equations, and some of theirs properties that will be useful in the following.

\subsection{The shallow water equations}

The 1D shallow water equations, introduced by Saint-Venant in \cite{SV1871} and derived in \cite{GP2001} from Navier-Stokes incompressible equations with a free moving boundary, describe the water propagation in a canal with rectangular cross-section and constant slope as follows:
\begin{equation}\label{eq:saint_venant}
\left\{\begin{array}{l}
\partial_t h + \partial_x (hv) = 0,\\
\\
\partial_t (hv) + \partial_x( h v^2 + \frac 12 g h^2) = g h (S_0-S_f),
\end{array}\right.
\end{equation}
with $h(x,t)$ the water height, $v(x,t)$ the water velocity at time $t$ and location $x$ along the canal, $g$ the gravity constant, $S_0$ the bed slope function and $S_f$ the friction slope function. The first equation comes from mass conservation and the second one from momentum balance.
For the purpose of this work, we assume a steady state friction on all canals and we assume horizontal canals with zero slope. Thus, the source term is zero.
 
We set $q=hv$ (the quantity $hv$ is often called the \textit{discharge} in shallow water theory, since it measures the flow rate of water past a point) and we reformulate  
system \eqref{eq:saint_venant} in vector form as 
\begin{equation}\label{eq:shallow_water_homog}
\partial_t U + \partial_x f(U) = 0,
\end{equation}
where
\begin{equation}\label{eq:shallow_water_flux}
U=\left(\begin{array}{c}
h \\ q
\end{array}\right), \quad 
f(U)= \left(\begin{array}{c}
hv \\ h v^2 + \frac 12 g h^2
\end{array}\right). 
\end{equation}
For smooth solutions, system \eqref{eq:shallow_water_homog} can equivalently be written in the quasilinear form
\begin{equation}\label{eq:sys_form}
\partial_t U + A(U)\partial_x U= 0,
\end{equation}
where the Jacobian matrix $A(U) = f^\prime(U)$ is
\begin{equation}
\label{eq:J_matrix}
A(U) = \left(\begin{array}{cc} 
0 & 1 \\
-v^2+g h & 2 v
\end{array}\right),
\end{equation}
with eigenvalues
\begin{equation}\label{eq:eigenvalues}
\lambda_1(U) = v-\sqrt{gh}, \quad \lambda_2(U) = v + \sqrt{gh}.
\end{equation}
Note that in general $\lambda_1$ and $\lambda_2$ can be of either sign.
When the velocity $v=q/h$ of the fluid is smaller than the speed $\sqrt{gh}$ of the gravity waves, that is $|v| < \sqrt{gh}$,
the flow is said to be \textit{fluvial} or \textit{subcritical} and then one has 
\begin{equation}\label{hyp:subcritical}
\lambda_1<0,\quad \lambda_2>0.
\end{equation}
Hence, under the subcritical condition \eqref{hyp:subcritical}, there are two waves propagating in opposite directions. 
The left and right characteristics are associated to $\lambda_1$ and $\lambda_2$ respectively. The ratio $ Fr = |v|/\sqrt{gh}$
is called the \textit{Froude number} and the flow is subcritical iff $\ds  Fr<1$.

%
\subsection{The standard Riemann problem for shallow-water equations.}
Here we are in particular interested in the solution of the Riemann problem:
\begin{equation}\label{eq:shallow_water_Rpb}
\left\{\begin{array}{l}
\partial_t U + \partial_x f(U) = 0,\\
\\
U(x,0) = \left\{
\begin{array}{ll}
U_l & \mbox{ if } x<0,\\ U_r & \mbox{ if } x>0,
\end{array}
\right.
\end{array}\right.
\end{equation}
where $U(x,0)=(h(x,0),q(x,0))$ is the initial condition and $U_l=(h_l,q_l)$ (resp. $U_r=(h_r,q_r)$) is the initial constant state on the left (resp. on the right) of the interface $x=0$.
The characteristic fields of the shallow water equations are genuinely nonlinear and so the Riemann problem always consists of two waves, each of which is either a shock or a rarefaction.
Under the \emph{subcritical flow condition} \eqref{hyp:subcritical}, there will be one left (with negative speed) and one right (with positive speed) going wave. In the sequel the left and right going waves are denoted by $l$-wave and $r$-wave, respectively. The solution to this Riemann problem consists of the $l$-wave and the $r$-wave separated by an intermediate state $\hat U=(\hat h,\hat q)$. We remark that the solution at the interface $x=0$ coincides with $\hat U$, which is the intersection point of the two functions $\phi_l$ and $\phi_r$ defined by 
\begin{equation}\label{eq:Rpb_1wave}
\phi_l(h;U_l) = \left\{
\begin{array}{ll}
v_l-2(\sqrt{gh}-\sqrt{g h_l}) & \mbox{ if } h<h_l \mbox{ (rarefaction)}\\
v_l - (h-h_l)\sqrt{g\frac{h+h_l}{2hh_l}} & \mbox{ if } h>h_l \mbox{ (shock wave)},
\end{array}
\right.
\end{equation}
and
\begin{equation}\label{eq:Rpb_2wave}
\phi_r(h;U_r) = \left\{
\begin{array}{ll}
v_r+2(\sqrt{gh}-\sqrt{g h_r}) & \mbox{ if } h<h_r \mbox{ (rarefaction)}\\
v_r + (h-h_r)\sqrt{g\frac{h+h_r}{2hh_r}} & \mbox{ if } h>h_r \mbox{ (shock wave)},
\end{array}
\right.
\end{equation}
which return the physically correct $\hat h$ and $\hat v$ intermediate values connecting the left  and right states with an entropic solution.

%
\section{Angle dependent conditions at the junction}\label{sec:Junction}
In this work, a junction is defined as the intersection of three channels. We assume that 1D shallow water equations hold on each canal of the network and we aim at deriving coupling conditions at the junction. These conditions enable to compute the intermediate states at the junction for the Riemann problem under consideration.
\subsection{Definition of the coupling conditions at the junction}\label{sec:Coupling_junction}
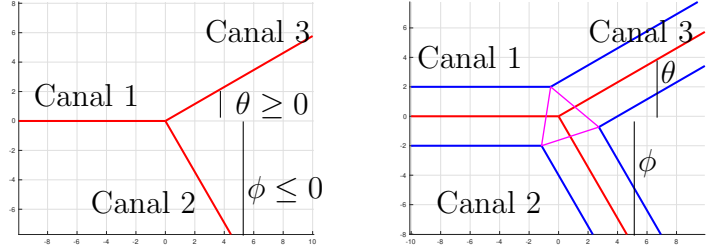
\begin{figure}[htbp!]\label{triangle}
\begin{overpic}
[width=0.4\textwidth]{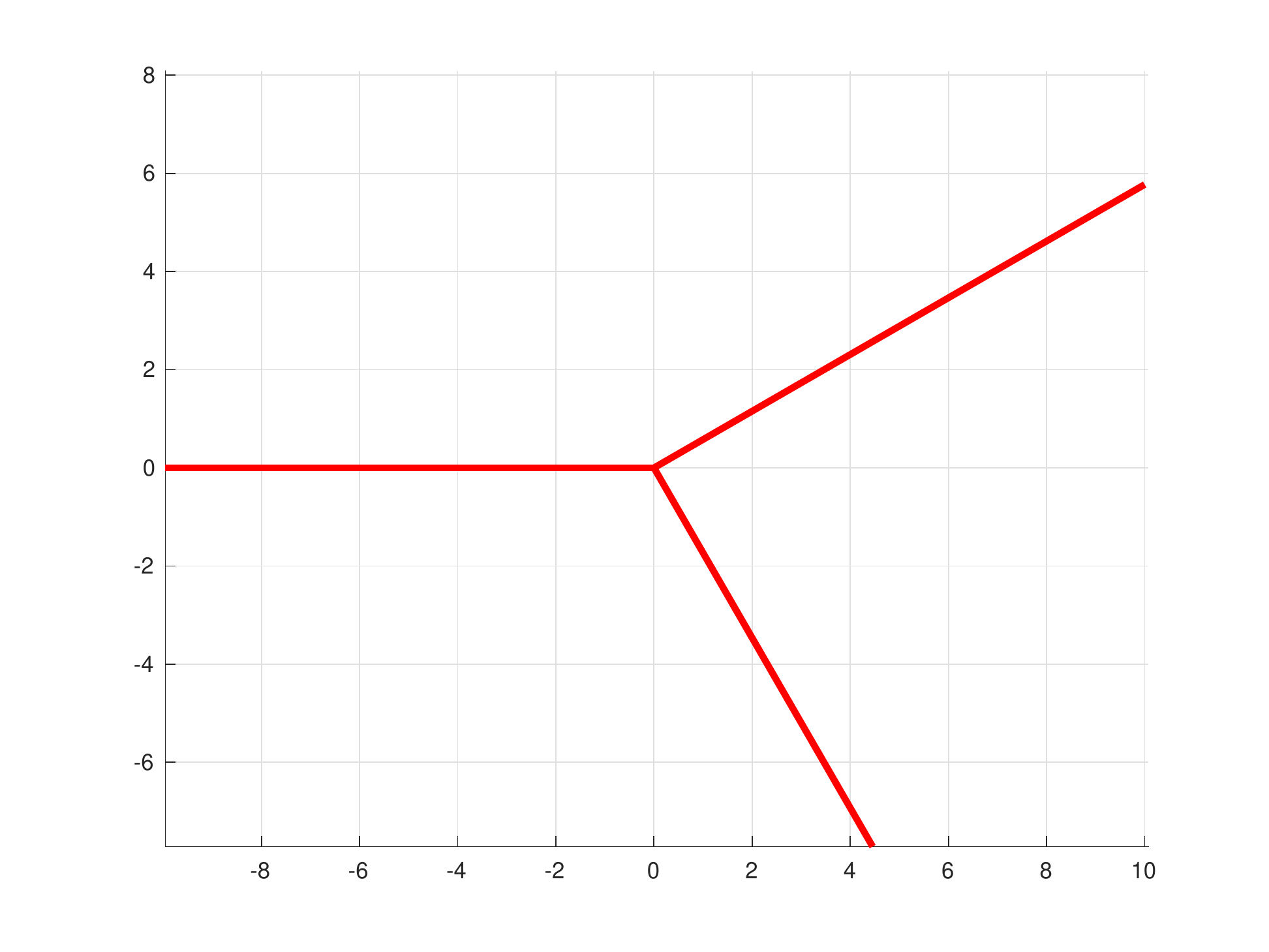}
\put(17,42){Canal 1}
\put(32,14){Canal 2}
\put(62,59){Canal 3}
\put(66,39){\line (0,1){7}}
\put(72,8){\line(0,1){30}}
\put(70,40){$\theta\geq 0$}
\put(73,18){$\phi\leq 0$}
\end{overpic} 
\begin{overpic}
[width=0.4\textwidth]{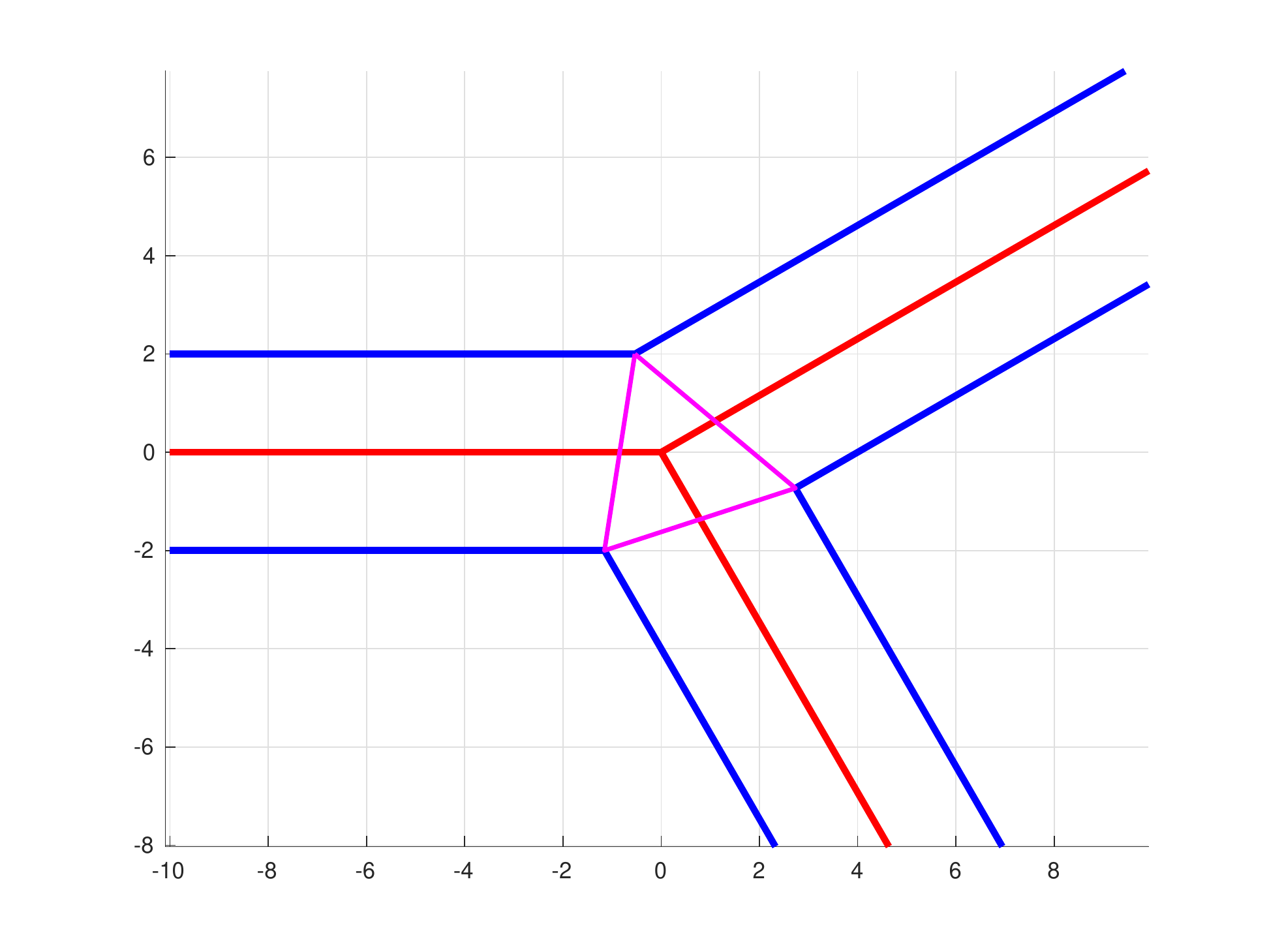}
\put(15,52){Canal 1}
\put(20,14){Canal 2}
\put(60,59){Canal 3}
\put(78,39){\line(0,1){15}}
\put(72,8){\line(0,1){30}}
\put(79,48){$\theta$}
\put(73,25){$\phi$}
\end{overpic} 
\caption{A 3 canal junction. On the left, the 1D set-up; on the right, the 2D  configuration. }\label{Fig:Geo}
\end{figure}

The channels will be labeled 1, 2 and 3 respectively, where channel 1 is assumed to be parallel to the $x$ axis. We fix the origin of the reference system in the point where the three channels intersect. Let $\theta$ and $\phi$ be the angles that channel 3 and 2 respectively form with the $x$ axis. We will assume that $\theta\geq 0$, while $\phi\leq 0$, obtaining the geometry in Fig. \ref{Fig:Geo} on the left. This is the one dimensional set up. 

Further, we will suppose that the channels can have different widths. Let then $2s_j, j=1, 2, 3$ be the  width of each channel. Therefore, we can think that the 1D setup is the core of a two dimensional junction, as shown in Fig. \ref{Fig:Geo} on the right.

Let $I_k$, $k=1,2,3$, be the interface separating the $k$-th channel from the junction. Let 
$U_k^*$, $k=1,2,3$,  denote the state variable in channel $k$ at the side of $I_k$ facing the channel obtained with the 1D solver used in the canal, while $U_k$, $k=1,2,3$, is the state variable at the side of $I_k$ facing the junction. The purpose of the junction Riemann solver is to compute $U_k$ given $U_k^*$. 
Since each state consists of the couple $(h,v)$, we need to find 6 unknowns at the junction. Three conditions are obtained finding the intermediate states of the one dimensional Riemann problem defined at each interface $I_k$ and in order to compute the three other missing data we shift to the 2D setting of Fig.\ref{Fig:Geo} on the right. We consider the triangle formed by the intersection points of the walls of the three channels, and to this two dimensional figure we apply conservation of mass and of the two components of momentum, which gives us the $3$ missing equations. 
Once the three states $U_k$, $k=1,2,3$, at the junction have been computed, we have at each interface $I_k$ the left and right states which are needed to compute the numerical flux at the boundary interfaces of the channels.

\subsubsection{Junction conditions coming from the Riemann solver.}

Let us begin with the $3$ equations coming from the Riemann solver.
We emphasize that, by convention,  the given configuration fixes channel 1 as entering the junction and channels 2 and 3 as leaving the junction.

Three  relations are obtained matching the unknowns $U_k$ at the junction with the data $U_k^*$ coming from the three channels through equations \eqref{eq:Rpb_1wave} and \eqref{eq:Rpb_2wave}. More precisely, 
\begin{equation}\label{eq:cond_riemann_pb}
\begin{array}{c}
v_1 = \phi_l(h_1;U^*_1)
\smallskip\\
v_2 = \phi_r(h_2;U^*_2)
\smallskip\\
v_3 = \phi_r(h_3;U^*_3).
\end{array}
\end{equation}
We note again that this construction requires a {\em fluvial regime}, in which only one wave exits the junction towards each of the three adjoining channels.

\subsubsection{Junction conditions coming from mass and momentum conservation.}
We now derive the  $3$ supplementary equations coming from conservation of mass and of the two components of momentum. For that purpose, we come back to the 2D configuration of the junction and we use the following notations
 \begin{itemize}
 \item $h$ denotes the height of water in the 2D configuration, 
\item  $\ds \mathbf{v}=(v_x,v_y)$, denotes the 2D velocity in the 2D junction domain, see Fig.\ref{Fig:Geo},
\item $\ds \mathbf{q}=h \mathbf{v}=h(v_x,v_y)$ denotes the 2D discharge,
\item $\ds \mathbf{q}_{k}=h_{k}(v_{x,k},v_{y,k})$, $k=1,2, 3$, denotes the average discharge on the edge of the junction triangle corresponding to channel $k$, see Fig.\ref{Fig:Geo}.
\end{itemize}
 We first recall the shallow-water equations in 2D, composed of the mass conservation equation and of the momentum conservation equation:
\begin{equation}\label{eq:saint_venant2D}
\left\{\begin{array}{l}
\partial_t h + \nabla \cdot (h\mathbf{v}) = 0,\\
\\
\partial_t (h\mathbf{v}) + \nabla \cdot( h \mathbf{v}\otimes \mathbf{v}) +\nabla( \ds\frac 12 g h^2) = 0.
\end{array}\right.
\end{equation}
In the following, we call $T$  the triangle formed by the intersection points of the walls of the three channels and its boundary  $\partial T$ is composed of three edges, denoted by $e_{k}$, $k=1,2,3$, see Fig.\ref{Fig:Geo}.
 
Mass conservation across the triangle $T$ with boundary $\partial T$ yields
\begin{equation}\label{eq:gen_mass}
\int_{\partial T}  \mathbf{q}  \cdot  \mathbf{n}=0, \text{ with }\mathbf{q}=h \mathbf{v},
\end{equation}
where $ \mathbf{n}$ is the outer normal of  $\partial T$, 
while the conservation of the two components of momentum gives the two relations
\begin{equation}\label{eq:gen_momentum1}
\int_{\partial T}  \left( v_{x}   \mathbf{q} + \frac 1 2 g h^2 \left(\begin{array}{c} 1 \\ 0 \end{array}\right) \right) \cdot  \mathbf{n}=0,
\end{equation}
and
\begin{equation}\label{eq:gen_momentum2}
\int_{\partial T}  \left( v_{y}   \mathbf{q} + \frac 1 2 g h^2\left(\begin{array}{c} 0 \\ 1 \end{array}\right) \right) \cdot  \mathbf{n}=0.
\end{equation}

Decomposing $\partial T$ as the sum of the three edges $e_{k}$, $k=1,2,3$, the three conditions at the junction, given by mass conservation and the two components of momentum conservation, can then be written as:
 \begin{subequations}
 \label{eq:2Dconditions}
\begin{align}
&  \sum_{k=1,2,3} \ell_{k}\mathbf{q}_{k}  \cdot  \mathbf{n}_{k}=0, \label{eq:2Dconditions:Mass}\\
& \sum_{k=1,2,3} \ell_{k}\left(v_{x,k}  \mathbf{q}_{k}   +\frac 1 2 g h_{k}^2
\left(\begin{array}{c} 1 \\ 0 \end{array}\right)
\right)\cdot  \mathbf{n}_{k}=0,  \label{eq:2Dconditions:Momentum1}\\
& \sum_{k=1,2,3} \ell_{k}\left(v_{y,k}  \mathbf{q}_{k}   +\frac 1 2 g h_{k}^2
\left(\begin{array}{c} 0 \\ 1 \end{array}\right)
\right)\cdot  \mathbf{n}_{k}=0,\label{eq:2Dconditions:Momentum2}
\end{align}
\end{subequations}
where $\ell_{k}$ is the length of the edge $e_{k}$  of the triangle,  $\mathbf{n}_{k}$ is the outer normal to  $e_{k}$ and $\mathbf{q}_{k}$ is the average of $\mathbf{q}$ on the side $e_{k}$ of the triangle.

\begin{figure}[htbp!]
\begin{overpic}
[width=0.9\textwidth]{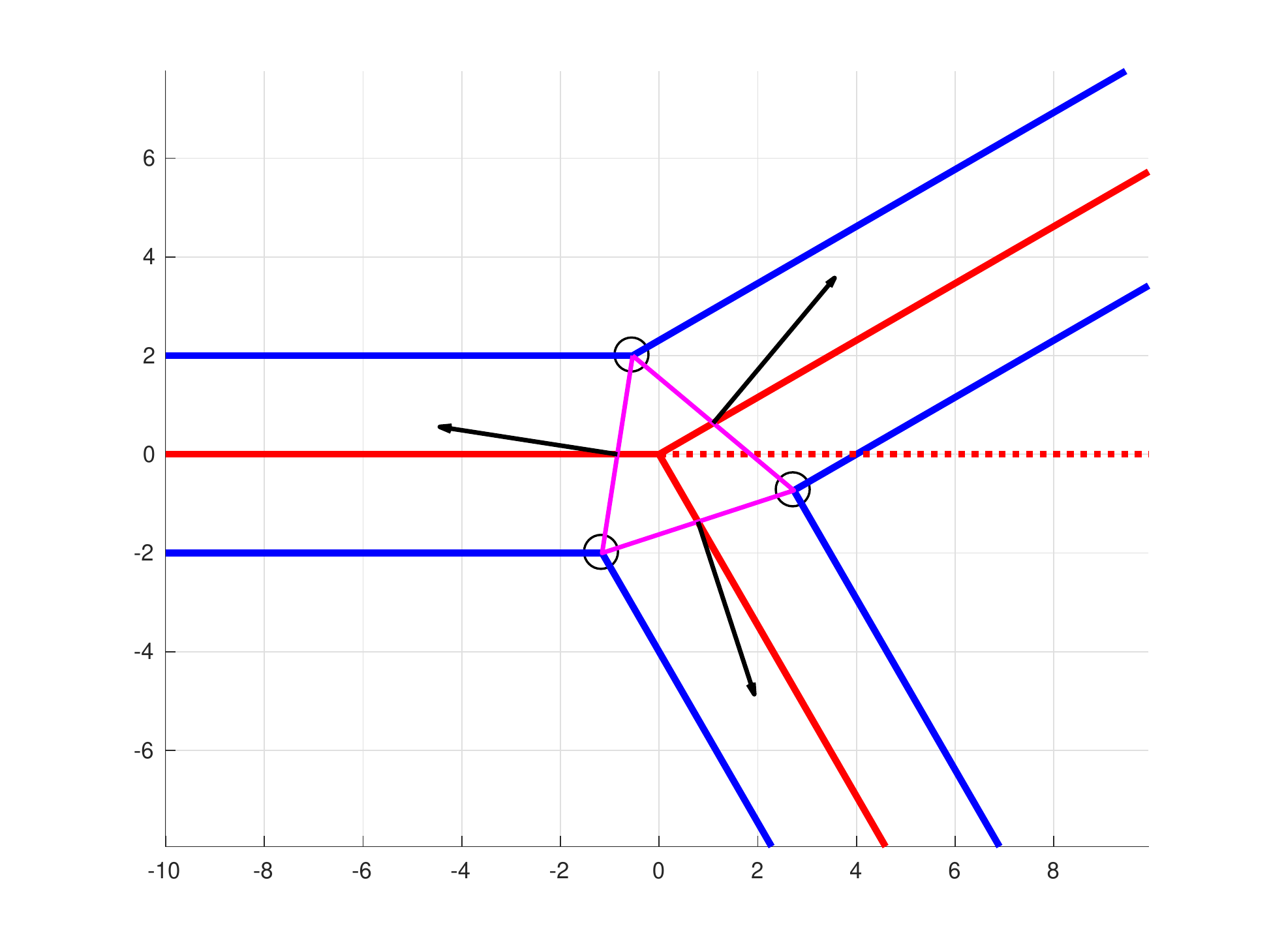}
\put(40,43){$\mathbf{n_1}$}
\put(53,25){$\mathbf{n_2}$}
\put(57,48){$\mathbf{n_3}$}
\put(20,42){Canal 1}
\put(58,14){Canal 2}
\put(75,59){Canal 3}
\put(78,39){\line(0,1){15}}
\put(72,8){\line(0,1){30}}
\put(79,50){$\theta\geq 0$}
\put(73,28){$\phi\leq 0$}
\put(42,26){$P_{12}$}
\put(46,53){$P_{13}$}
\put(66,35){$P_{23}$}
\end{overpic} 
\caption{A 3 canal junction. Illustration of the geometrical notations. Parameters are $\ds s_{1}=s_{2}=s_{3}=2$,  $\theta=\frac{\pi}{6}$ and $\phi=-\frac{\pi}{3}$}\label{fig:3canals-junction}
\end{figure}
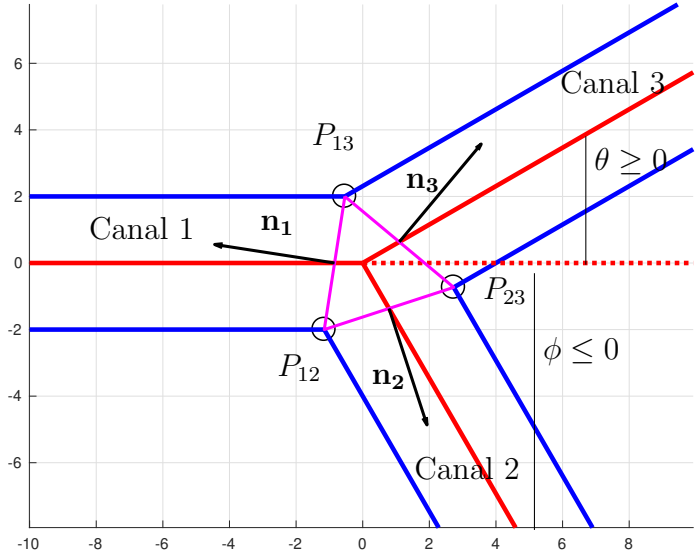

To specify all quantities appearing in  system \eqref{eq:2Dconditions}, we need to compute  the normals  $\mathbf{n}_{k}$   to the sides of the triangle and their lengths $\ell_{k}$. To fix notation, refer to Fig. \ref{fig:3canals-junction}.

To begin with, we need to give the coordinates of the intersection points of the walls, namely points $P_{12}$, $P_{13}$ and $P_{23}$ that are displayed on Fig. \ref{fig:3canals-junction}. Let us recall that $\theta$ and $\phi$ are the angles of canals 2 and 3 with the $x$-axis, while $2 s_{k}$ is the section of canal $k$.

 The equations for the straight lines composing the 1D skeleton of Fig.\ref{Fig:Geo} written in parametric form are
\begin{align*}
y_1 & = t_1[1,0]^T, \\
y_2 & = t_2[\cos(\phi),\sin(\phi)]^T, \\
y_3 & = t_3[\cos(\theta),\sin(\theta)]^T,
\end{align*}
with $t_k\in \R, k=1, 2, 3$. Then, to obtain the walls of the channels, i.e, to construct the 2D setting of Fig. \ref{fig:3canals-junction}, we just need to write the equations of the two straight lines parallel to the axis $y_k$ at the center of the channel, and at a distance $\pm s_k$ from the axis, for each channel. The walls of the three channels are
\begin{align}
y_1^\pm & = t_1[1,0]^T \pm s_1[0,1]^T, \nonumber \\
y_2^\pm & = t_2[\cos(\phi),\sin(\phi)]^T \pm s_2[-\sin(\phi),\cos(\phi)]^T, \\
y_3^\pm & = t_3[\cos(\theta),\sin(\theta)]^T \pm s_3[-\sin(\theta),\cos(\theta)]^T. \nonumber
\end{align}
The triangle in Fig. \ref{fig:3canals-junction} across which the 2D interaction occurs is obtained intersecting the straight lines defining the walls of the channels. More precisely, $P_{13}$, is the intersection of $y_1^+$ with $y_3^+$, $P_{12}$ is defined by the intersection of $y_1^-$ with $y_2^-$, and the last point $P_{23}$ lies at the intersection of $y_3^-$ and $y_2^+$. 
We obtain, 
\begin{equation}\label{e:P13}
P_{13}= \left( \frac{s_1 \cos \theta- s_3}{\sin \theta}, s_1 \right), \quad \theta\ne 0.
\end{equation}
If $\theta = 0$, the system has a solution only provided $s_1=s_3$, and the two straight lines actually coincide. In this case we define $P_{13}=(0,s_1)$. 

Analogously,
\begin{equation}\label{e:P12}
P_{12}= \left( \frac{-s_1\cos \phi+s_2}{\sin \phi}, -s_1 \right), \quad \phi\ne 0.
\end{equation}
If $\phi=0$, we must have $s_1=s_2$, and we pick $P_{12}=(0,-s_1)$. With this approach, we cannot treat the case in which both $\phi=\theta=0$, unless we consider the two channels $y_2$ and $y_3$ superposed one on top of the other. We will see in the next subsection how to extend the construction also to the case $\phi=\theta=0$.

Finally, 
\begin{equation}\label{e:P23}
P_{23}= \left( \frac{s_3\cos \phi+s_2\cos \theta}{\sin (\theta-\phi)}, \frac{s_3\sin \phi+s_2 \sin \theta}{\sin (\theta-\phi)} \right).
\end{equation}
The quantity $\sin(\phi-\theta)$ can be zero either for $\phi=\theta=0$, in which case the two pipes coincide, or when $\phi=-\pi/2$ and $\theta=\pi/2$. Now you have solutions only for $s_3=s_2$, which means that $y_3^-$ and $y_2^+$ coincide, and we fix the intersection point to $P_{23}=(s_3,0)$.

We will analysize, and extend, the particular cases $\theta=0, \phi=0, (\theta,\phi)=(\pi/2,-\pi/2)$ in the following section. 

Once the points $P_{13}, P_{12}, P_{23}$ are defined, we can compute all quantities  $\mathbf{n}_{1}$,  $\mathbf{n}_{2}$,  $\mathbf{n}_{3}$, $\ell_{1}$, $\ell_{2}$ and $\ell_{3}$ depending on the geometry appearing in 
\eqref{eq:2Dconditions}.
The length of the sides is
\begin{equation}
\ell_1=||P_{13}-P_{12}||, \qquad \ell_2=||P_{23}-P_{12}||, \qquad \ell_3=||P_{23}-P_{13}||,
\end{equation}
and the normals are
\begin{equation*}
 \mathbf{n}_{1}=\frac{1}{\ell_{1}}
 \left(
\begin{array}{c}
  -2s_{1}
\medskip
   \\
\ds\frac{s_{1} \sin (\theta+\phi)-s_{2}\sin \theta -s_{3} \sin\phi}{\sin \phi\sin \theta} 
\end{array}
\right), 
\end{equation*}
\begin{equation*}
 \mathbf{n}_{2}=\frac{1}{\ell_{2}}
 \left(
\begin{array}{c}
  s_{1}+\ds\frac{s_{2}\sin \theta +s_{3} \sin\phi}{\sin (\theta-\phi)} 
\medskip
   \\
\ds\ds-\frac{s_{2}\cos \theta +s_{3} \cos\phi}{\sin (\theta-\phi)} -
\frac{s_{1} \cos \phi-s_{2}}{\sin \phi} 
\end{array}
\right), 
\end{equation*}
\begin{equation*}
 \mathbf{n}_{3}=\frac{1}{\ell_{3}}
 \left(
\begin{array}{c}
   s_{1}-\ds\frac{s_{2}\sin \theta +s_{3} \sin\phi}{\sin (\theta-\phi)} 
\medskip
   \\
   \ds\frac{s_{2}\cos \theta +s_{3} \cos\phi}{\sin (\theta-\phi)} -\frac{s_{1}\cos \theta-s_{3}}{\sin \theta}
\end{array}
\right).
\end{equation*}

\begin{remark}\label{alignment} The construction is well defined as long as the triangle formed by $P_{13}, P_{12}, P_{23}$ is non degenerate. We say that the triangle is degenerate when the three points lie on the same straight line. 
Straightforward computations show 
 that this occurs when $\ds \det ( \mathbf{n}_{1}, \mathbf{n}_{3})=0$
 which is equivalent to the particular combination 
\begin{equation}\label{eq:EqMagali}
 (s_{1}\sin( \theta-\phi) +s_{3}\sin(\phi) - s_{2}\sin(\theta))^2+4 s_{2}s_{3} \sin(\phi)\sin(\theta) = 0.
\end{equation}
\end{remark}

In the frame of reference we have chosen, the discharge in the three canals can be written as
\begin{equation}
\begin{array}{ccc}
\mathbf{q}_1 = q_1 \left(\begin{array}{c} 1 \\ 0\end{array}\right), & \mathbf{q}_2 =q_2  \left(\begin{array}{c} \cos\phi \\ \sin\phi\end{array}\right), &\mathbf{q}_3 =q_3 \left(\begin{array}{c} \cos\theta \\ \sin\theta\end{array}\right),
\end{array}
\end{equation}
where $q_k = \|\mathbf{q}_k\|$.

Let $v_k=q_k/h_k$ be the velocity along the $k-$th channel.
 Then the conservation laws \eqref{eq:2Dconditions} across the junction can be written as
 \begin{subequations}\label{eq:Conservation}
  \begin{align}
& \ell_1 h_{1}v_1  \left(\begin{array}{c} 1 \\ 0\end{array}\right)\cdot \mathbf{n}_{1} + \ell_2 h_{2}v_2  \left(\begin{array}{c} \cos\phi \\ \sin\phi\end{array}\right) \cdot \mathbf{n}_{2} + \ell_3 h_{3} v_3 \left(\begin{array}{c} \cos\theta \\ \sin\theta\end{array}\right)\cdot \mathbf{n}_{3} = 0,  \label{eq:Conservation:Mass}\\
&  \ell_{1}\left(h_1  v_1^2 \left(\begin{array}{c} 1 \\ 0\end{array}\right)  +\frac 1 2 g h_{1}^2
\left(\begin{array}{c} 1 \\ 0 \end{array}\right)
\right)\cdot  \mathbf{n}_{1} +
 \ell_{2}\left(h_2 v_2^2 \cos\phi  \left(\begin{array}{c} \cos\phi \\ \sin\phi\end{array}\right)   +\frac 1 2 g h_{2}^2
\left(\begin{array}{c} 1 \\ 0 \end{array}\right)
\right)\cdot  \mathbf{n}_{2},\label{eq:Conservation:Momentum1}  \\
& \qquad \qquad \qquad + \ell_{3}\left(h_3 v^2_3 \cos \theta \left(\begin{array}{c} \cos\theta \\ \sin\theta\end{array}\right)   +\frac 1 2 g h_{3}^2
\left(\begin{array}{c} 1 \\ 0 \end{array}\right)
\right)\cdot  \mathbf{n}_{3}=0, \notag \\
&  \ell_{1}\left(
\frac 1 2 g h_{1}^2
\left(\begin{array}{c} 0 \\ 1 \end{array}\right)
\right)\cdot  \mathbf{n}_{1} +
 \ell_{2}\left(h_2 v_2^2 \sin \phi  \left(\begin{array}{c} \cos\phi \\ \sin\phi\end{array}\right)   +\frac 1 2 g h_{2}^2
\left(\begin{array}{c} 0 \\ 1 \end{array}\right)
\right)\cdot  \mathbf{n}_{2}\label{eq:Conservation:Momentum2} \\
& \qquad \qquad \qquad  +\ell_{3}\left(h_3 v_3^2 \sin \theta \left(\begin{array}{c} \cos\theta \\ \sin\theta\end{array}\right)   +\frac 1 2 g h_{3}^2
\left(\begin{array}{c} 0 \\ 1 \end{array}\right)
\right)\cdot  \mathbf{n}_{3}=0,\notag
\end{align}
\end{subequations}
where we used the fact that the axis of channel 1 is parallel to the $x$ axis.

\subsubsection{Solutions for the whole system of equations at the junction.}

Combining the three  equations \eqref{eq:Conservation} with the three equations \eqref{eq:cond_riemann_pb}
coupling the states in the junction with the 1D channels, we find a system of 6 non linear equations at the junction, whose solution is given by the three intermediate states $U_k=(h_{k}, v_{k})$, $k=1,2,3$.

\begin{remark}
If we consider  a stationary solution of \eqref{eq:saint_venant}  such that the velocity is null and the height is constant in space, i.e. $h_k^*=\bar h$ and $v_k^*=0$, $k=1,2,3$ then $h_k=\bar h$ and $v_k=0$ $k=1,2,3$ is a trivial solution  of system \eqref{eq:Conservation}-\eqref{eq:cond_riemann_pb} since $\ds \mathbf{n}_{1}+\mathbf{n}_{2}+\mathbf{n}_{3} =0$. This means that the coupling condition at the junction preserves the lake at rest stationary solution on the whole network.
\end{remark}

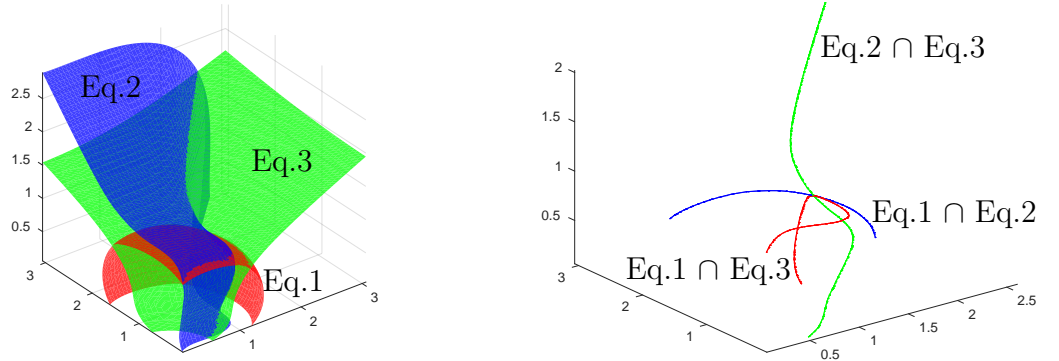
\begin{figure}[htbp!]
 \subfloat[Zero surfaces  of the three equations in \eqref{eq:Conservation}.]{
\begin{overpic}
[width=0.6\textwidth]{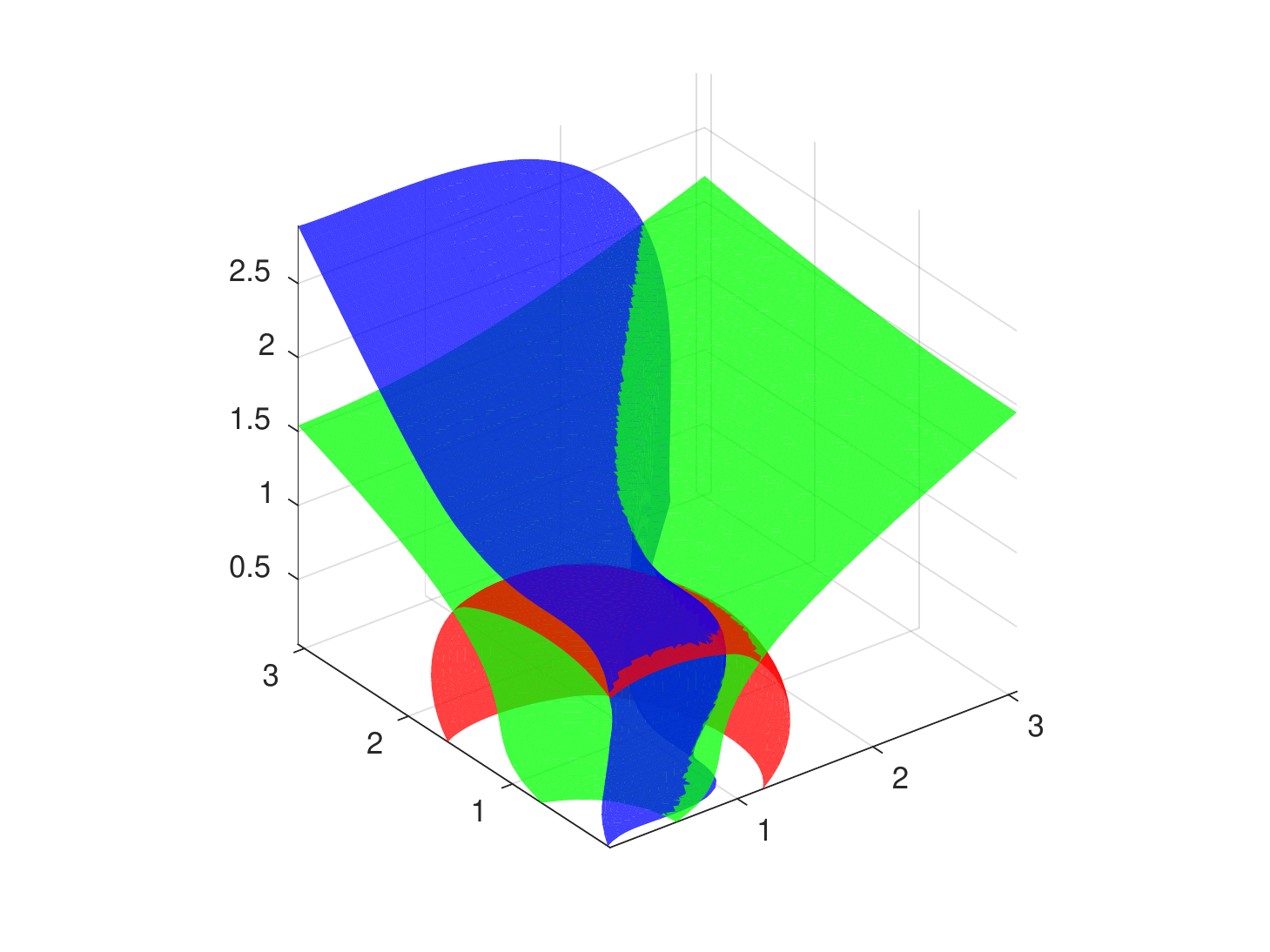}
\put(62,19){Eq.1}
\put(30,53){Eq.2}
\put(60,40){Eq.3}
\end{overpic} }
\subfloat[Intersection point of the three zero surfaces.]{
\begin{overpic}
[width=0.6\textwidth]{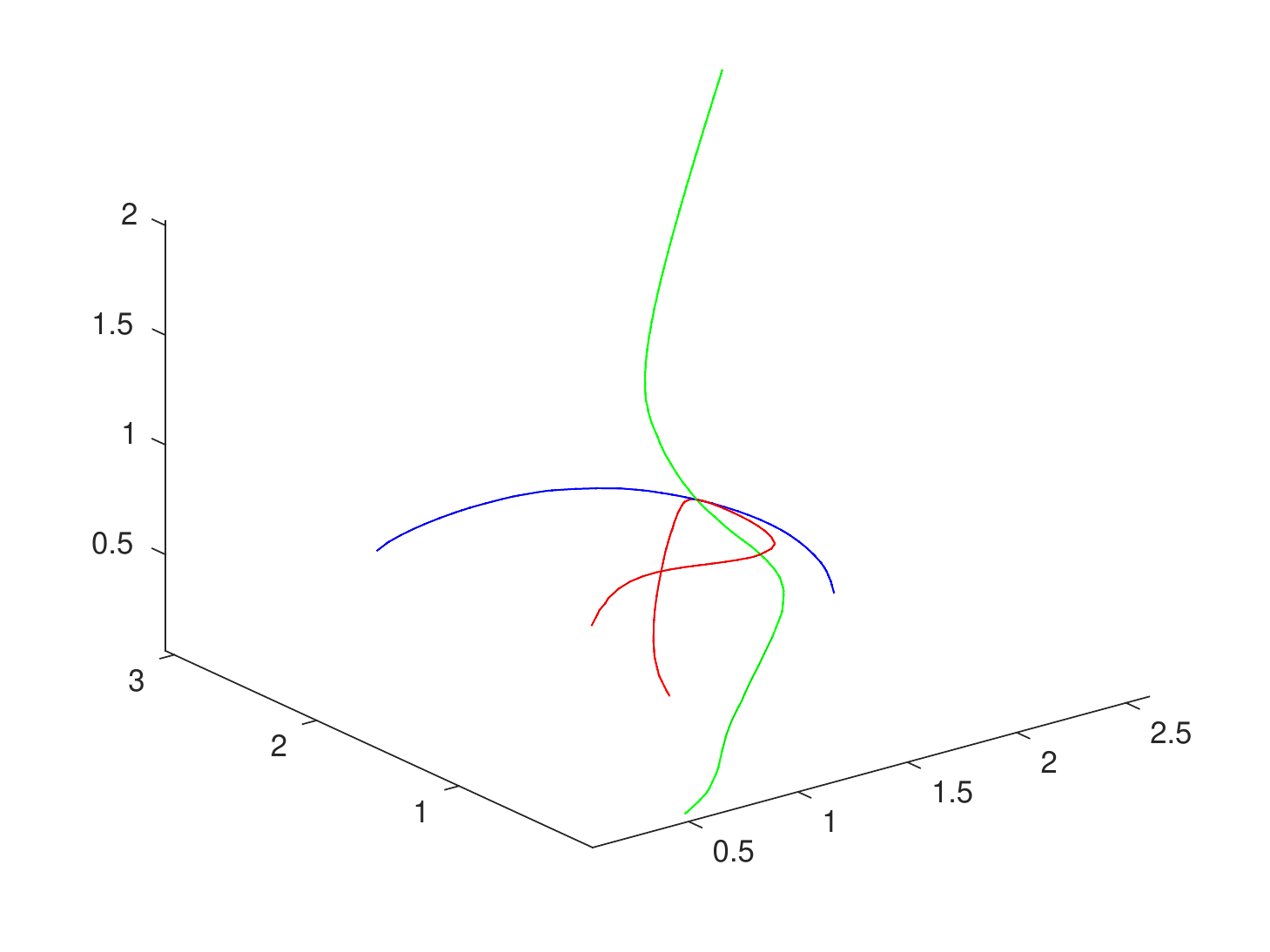}
\put(65,31){Eq.1 $\cap$ Eq.2}
\put(22,21){Eq.1 $\cap$ Eq.3}
\put(56,60){Eq.2 $\cap$ Eq.3}
\end{overpic} }
\caption{Graphic illustration of the existence of a unique solution of system  \eqref{eq:cond_riemann_pb}-\eqref{eq:Conservation} for the parameters $s_1=s_2=s_3=1$, $\theta=-\phi=\pi/6$ and $U^*_1=(1.5,0)$, $U^*_2=U^*_3=(1,0)$.}\label{fig:SurfSol}
\end{figure}

Substituting $v_1, v_2, v_3$ from \eqref{eq:cond_riemann_pb} into \eqref{eq:Conservation}, we find a system of three non linear equations in the three unknowns $h_1, h_2$ and $h_3$ which gives the solution at the junction. Once the parameters  $s_1, s_2, s_3$, $\theta,\phi$ and $U^*_1, U^*_2, U^*_3$ are fixed, these three equations define three hypersurfaces whose zeros surfaces can be plot in the plane $h_1, h_2$ and $h_3$ (Figure \ref{fig:SurfSol}). The intersection of this surfaces is the required solution. 
In Figure \ref{fig:SurfSol}-right we
show an example.

\subsection{Special cases and extensions}\label{sec:j_special}
In this section, we consider three particular cases.

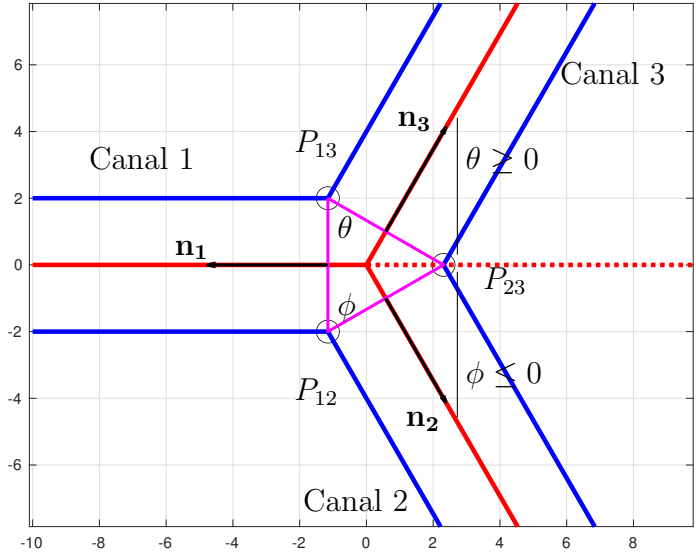
\begin{figure}[htbp!]\label{fig:3canals-junction-simple}
\begin{overpic}
[width=0.9\textwidth]{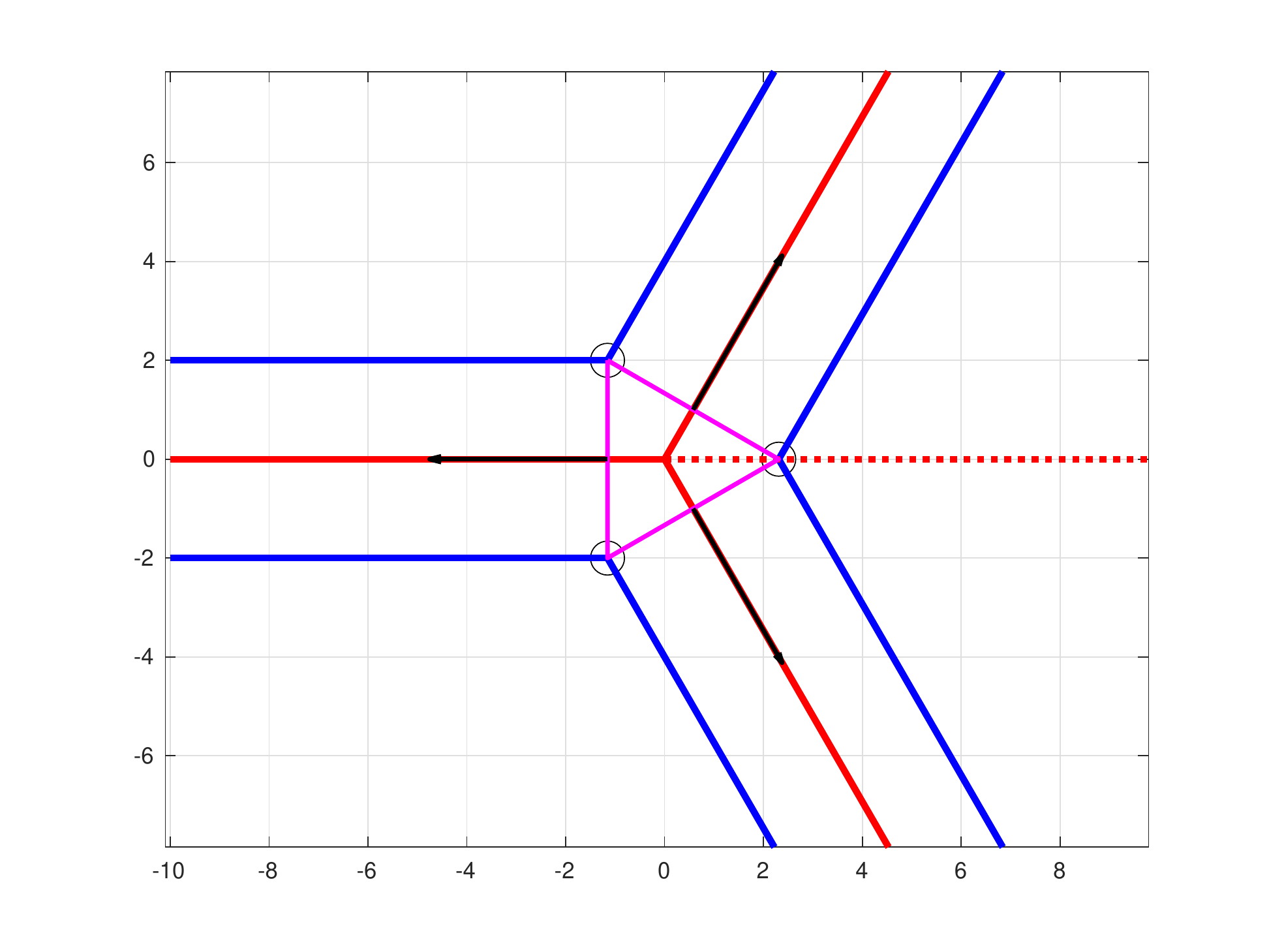}
\put(30,40){$\mathbf{n_1}$}
\put(57,20){$\mathbf{n_2}$}
\put(56,55){$\mathbf{n_3}$}
\put(20,50){Canal 1}
\put(45,10){Canal 2}
\put(75,60){Canal 3}
\put(63,40){\line(0,1){16}}
\put(63,21){\line(0,1){17}}
\put(64,50){$\theta\geq 0$}
\put(49,42){$\theta$}
\put(49,33){$\phi$}
\put(64,25){$\phi\leq 0$}
\put(44,23){$P_{12}$}
\put(44,52){$P_{13}$}
\put(66,36){$P_{23}$}
\end{overpic} 
\caption{A particular case : junction where the directions of the channels are  perpendicular to the sides of the triangle. Parameters are $\ds s_{1}=s_{2}=s_{3}=2$ and $\theta=-\phi=\frac{\pi}{3}$}\label{Tperp}
\end{figure}

We start with the simplified case in which the canals are orthogonal to the sides of the triangle. In that case, the junction is defined uniquely by the three sections and the equations \eqref{eq:Conservation} simplify loosing the dependency on the angles.
In this case, see Fig. \ref{Tperp}, it is easy to see that the angles $\theta$ and $\phi$ defining the skeleton of the junction coincide with the angles labelled $\theta$ and $\phi$ internal to the triangle in Fig. \ref{Tperp}, and the length of the sides coincides with the amplitude of the channels, namely $l_k=2s_k$, $k=1,2, 3$.

Then, it is straightforward  to see that the sections depend on the angles through the following relations
\begin{equation}\label{eq:OrthoPipesCondition}
\begin{array}{c}
s_2 \sin \phi + s_3 \sin \theta = 0,
\smallskip\\
s_1=s_2\cos \phi+s_3 \cos \theta.
\end{array}
\end{equation}
Since in the present case, $\mathbf{q}_k$ is parallel to $\mathbf{n}_k$, 
equation \eqref{eq:2Dconditions:Mass} becomes 
\begin{equation}\label{moment}
- s_{1}  q_{1}+s_{2}  q_{2}+s_{3}  q_{3}=0.
\end{equation}
Equation \eqref{eq:2Dconditions:Momentum1}-\eqref{eq:2Dconditions:Momentum2}, corresponding to the conservation of momentum at the junction in 2D give:
  \begin{equation}\label{mom1}
  \left( \frac{q_{1}^2}{h_{1}}+\frac 1 2 g h_{1}^2\right) s_{1}=\left( \frac{q_{2}^2}{h_{2}}+\frac 1 2 g h_{2}^2\right)  s_{2} \cos \phi+\left( \frac{q_{3}^2}{h_{3}}+\frac 1 2 g h_{3}^2\right)  s_{3}  \cos \theta
  \end{equation}
  and 
  \begin{equation}\label{mom2}
 0=\left( \frac{q_{2}^2}{h_{2}}+\frac 1 2 g h_{2}^2\right)  s_{2} \sin \phi+\left( \frac{q_{3}^2}{h_{3}}+\frac 1 2 g h_{3}^2\right)  s_{3}  \sin \theta.
  \end{equation}
  
  Using the identities in \eqref{eq:OrthoPipesCondition}, we can rewrite \eqref{mom1} and \eqref{mom2} as
  $$  \frac{q_{1}^2}{h_{1}}+\frac 1 2 g h_{1}^2= \frac{q_{2}^2}{h_{2}}+\frac 1 2 g h_{2}^2= \frac{q_{3}^2}{h_{3}}+\frac 1 2 g h_{3}^2.$$
  
  Therefore, conservation of mass and of the two components of momentum at the junction in this particular case yield
\begin{equation}\label{eq:conditions_simplecase}
\left\{\begin{array}{l}
-s_{1}  q_{1}+s_{2}  q_{2}+s_{3}  q_{3}=0,\\
\\
\ds \frac{q_{1}^2}{h_{1}}+\frac 1 2 g h_{1}^2=\ds \frac{q_{2}^2}{h_{2}}+\ds\frac 1 2 g h_{2}^2= \frac{q_{3}^2}{h_{3}}+\frac 1 2 g h_{3}^2.
\end{array}\right.
\end{equation}  
Note that in this case, the junction conditions do not depend on the angles with which the canals intersect.
We thus recover the equal energy condition at the junction used by several authors, see \cite{HS2013} and references there in. This condition derives from the 2D momentum conservation at the junction, but we stress that it holds only for the particular case in which the channels are orthogonal to the sides of the triangle defining the junction. 

Moreover, tedious but straightforward calculations show that equations \eqref{eq:OrthoPipesCondition} imply the geometry in Fig. \ref{Tperp}. Since the conservation condition can be multiplied by a constant without changing the result, we see that for each pair of angles $\phi$ and $\theta$ there exists a one parameter set of  sections $\lambda (s_{1}, s_2, s_{3}),$ with $\lambda>0$  for which momentum conservation coincides with energy conservation. 
\begin{remark}
The derivation leading to \eqref{eq:conditions_simplecase} proves that the present discussion is actually an {\em extension} of the junction conditions based on energy conservation. Only in the case of the particular combination of parameters satisfying \eqref{eq:OrthoPipesCondition}, the junction Riemann solver does {\em not depend} on the angles between the pipes. In particular, if all sections are equal, \eqref{eq:OrthoPipesCondition} implies that conservation of momentum coincides with conservation of energy only in the case $\theta=\ds\frac{\pi}{3}=-\phi$ see Fig~\ref{Tperp}.
\end{remark}

\medskip

We now consider the cases $\theta=\phi=0$ and $\theta=-\phi=\pi/2$ which where excluded in the generic case described in section \ref{sec:Junction}. We call~:
\begin{itemize}
\item T-junction: $\theta=-\phi=\pi/2$. 
\item Straight channel: $\theta=\phi=0$.
\end{itemize}

\subsubsection{T-junction}
\begin{figure}[htbp!]
 \subfloat[Case when $s_{2}=s_{3}$\label{fig:Tegal}]{
\begin{overpic}
[width=0.5\textwidth]{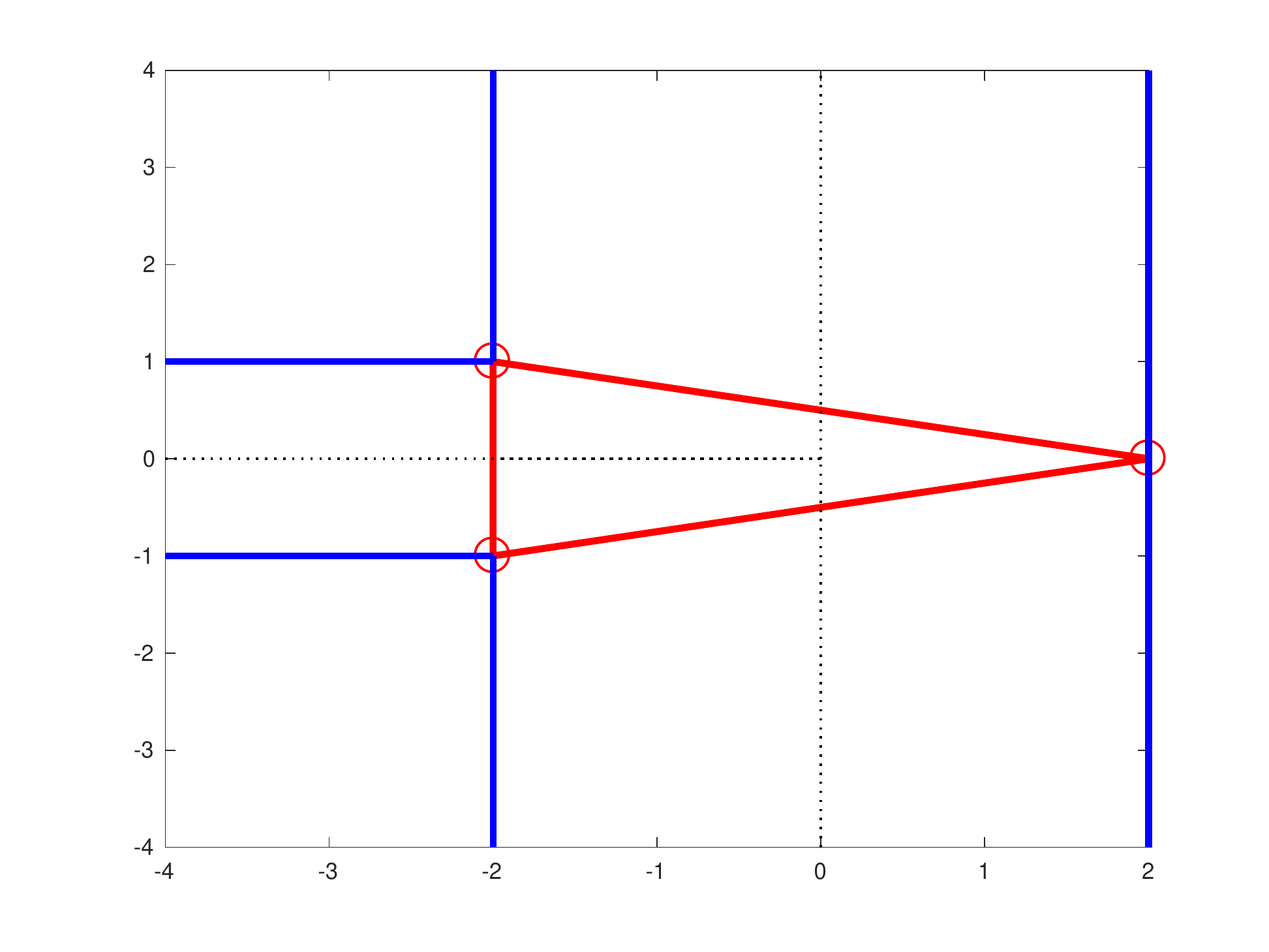}
\put(15,41){Canal 1}
\put(67,15){Canal 2}
\put(67,60){Canal 3}
\put(27,24){$P_{12}$}
\put(28,52){$P_{13}$}
\put(80,42){$P_{23}$}
\put(20,32){\line(0,1){6}}
\put(23,35){$\displaystyle s_{1}$}
\put(40,65){\line(1,0){24}}
\put(48,60){$\displaystyle s_{3}$}
\put(40,10){\line(1,0){24}}
\put(52,13){$\displaystyle s_{2}$}
\end{overpic} }
\subfloat[Case when $s_{2}\neq s_{3}$\label{fig:Tdiff}]{
\begin{overpic}
[width=0.5\textwidth]{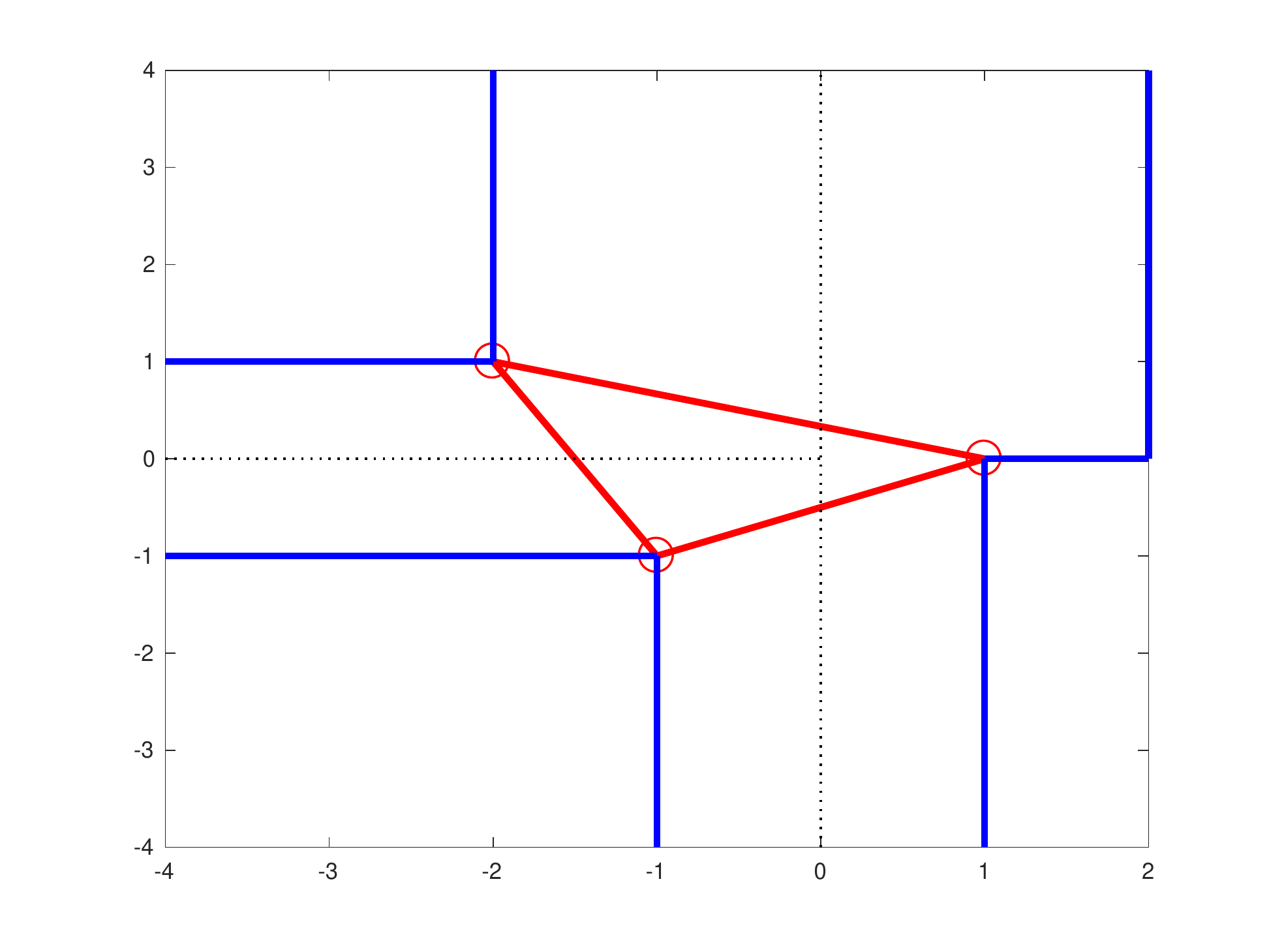}
\put(15,41){Canal 1}
\put(67,15){Canal 2}
\put(67,60){Canal 3}
\put(40,24){$P_{12}$}
\put(28,52){$P_{13}$}
\put(80,42){$P_{23}$}
\put(20,32){\line(0,1){6}}
\put(23,35){$\displaystyle s_{1}$}
\put(40,65){\line(1,0){24}}
\put(48,60){$\displaystyle s_{3}$}
\put(53,10){\line(1,0){11}}
\put(55,13){$\displaystyle s_{2}$}
\end{overpic} }
\caption{Junction - Particular case  of the T-junction. Parameters are $\theta=-\phi=\frac{\pi}{2}$ and $s_{1}=1$, $s_{2}=s_{3}=2$ (on the left) and $s_{1}=s_{2}=1$, $s_{3}=2$ (on the right).}\label{fig:Tjunction_cp}
\end{figure}
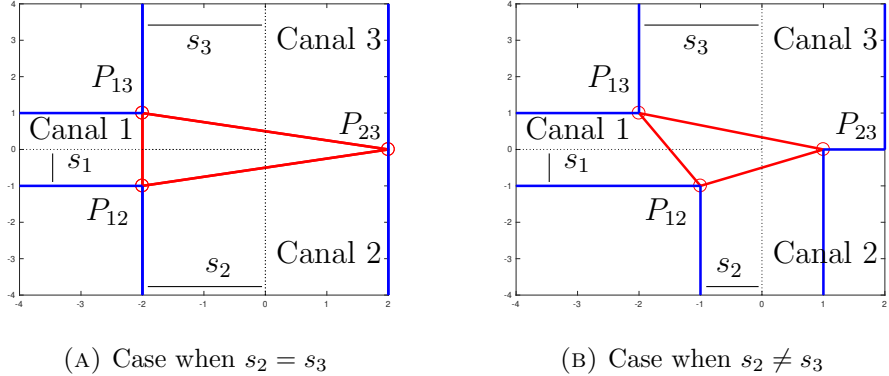

In the case of  a T-junction, for which the angles are equal to  $\theta=-\phi=\pi/2$ and $s_{2}=s_{3}$, the points $P_{12}$, $P_{13}$ and 
$P_{23}$ can be defined as in Fig.\ref{fig:Tjunction_cp}-(A), namely
\begin{equation*}
P_{12}=\left(
\begin{array}{c}
-s_{2}
\smallskip
 \\
-s_{1}
\end{array}
\right), 
\,
P_{13}=\left(
\begin{array}{c}
-s_{2}
\smallskip
   \\
  s_{1}
\end{array}
\right), 
\,
P_{23}=\left(
\begin{array}{c}
s_{2}
\smallskip
   \\
0
\end{array}
\right).
\end{equation*}
Equations \eqref{eq:Conservation} reduce to
\begin{equation*}
\left\{
\begin{aligned}
& -s_{1} h_1 v_{1}  + s_2 h_2 v_{2} + s_2 h_3 v_{3}  = 0,
\medskip\\
&  - 2\left( h_{1}v_{1}^2    +\frac 1 2 g h_{1}^2 \right)
 +
 \frac 1 2  g h_{2}^2
 + \frac 1 2  g h_{3}^2
=0, \\
& -\left(h_{2}v_{2}^2    +\frac 1 2 g h_{2}^2
\right)
+ h_{3}v_{3}^2   +\frac 1 2 g h_{3}^2
=0.
\end{aligned}
\right.
\end{equation*}
When $s_{2} \neq s_{3}$, the point $P_{23}$ can be defined as in Fig.\ref{fig:Tjunction_cp}-(B), then 
\begin{equation*}
P_{12}\left(
\begin{array}{c}
-s_{2}
\smallskip
 \\
-s_{1}
\end{array}
\right), 
\,
P_{13}\left(
\begin{array}{c}
-s_{3}
\smallskip
   \\
  s_{1}
\end{array}
\right), 
\,
P_{23}\left(
\begin{array}{c}
\min(s_{2},s_{3})
\smallskip
   \\
0
\end{array}
\right).
\end{equation*}

\subsubsection{Straight channel}
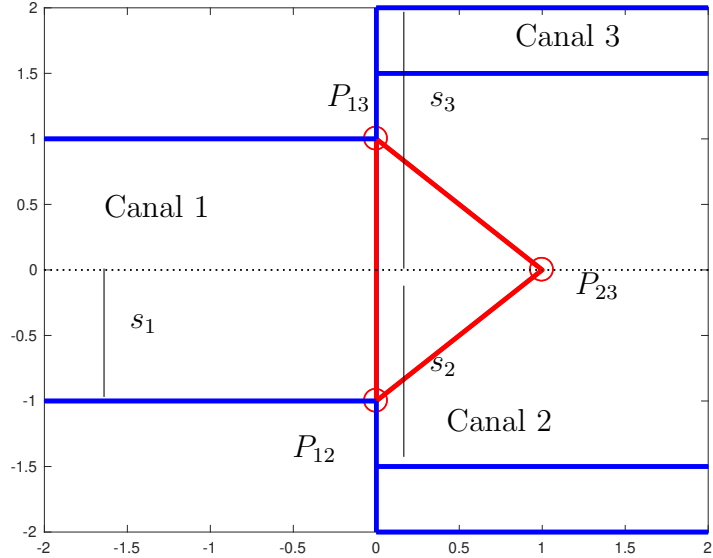
\begin{figure}[htbp!]
\begin{overpic}
[width=0.9\textwidth]{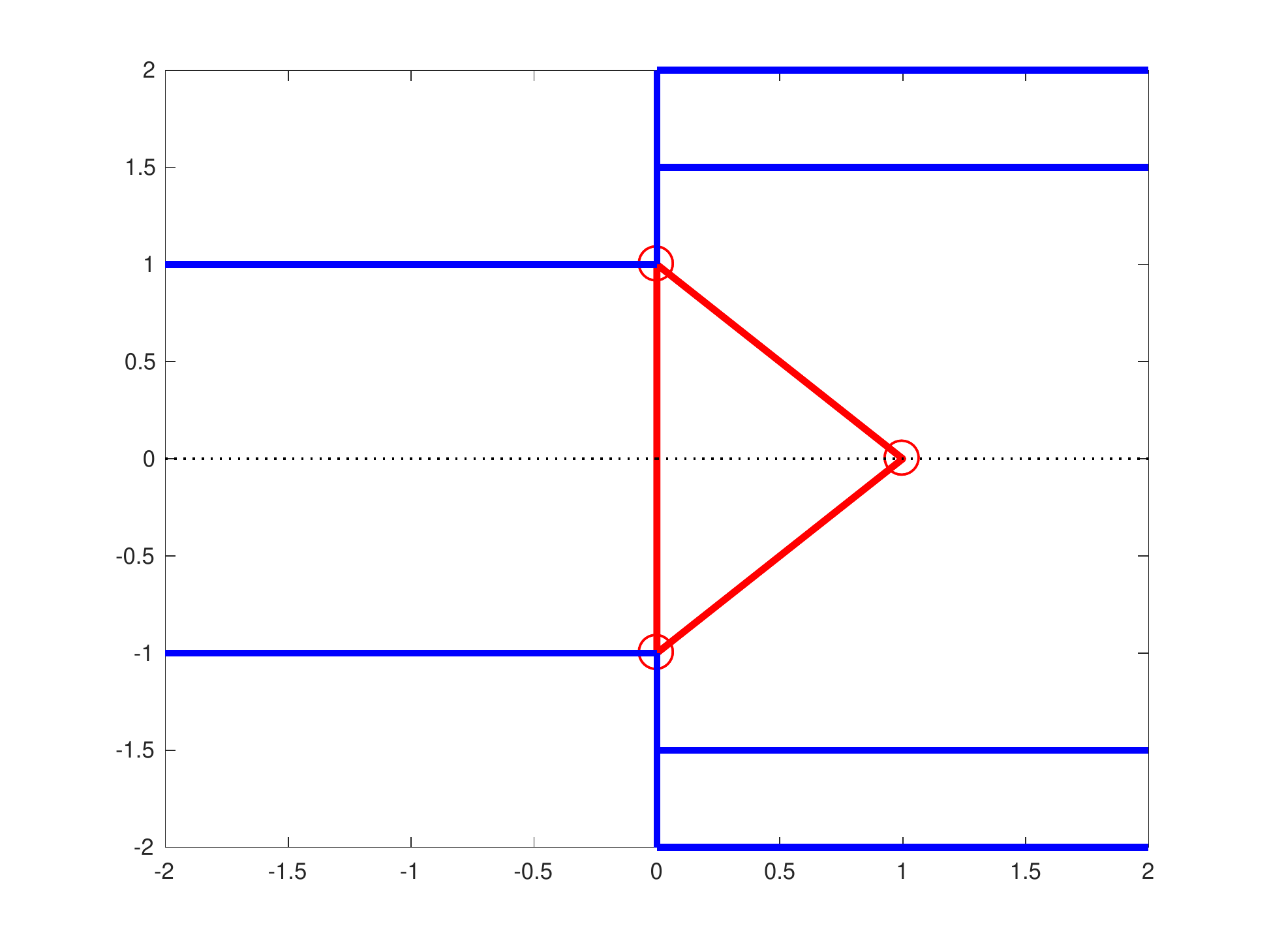}
\put(20,45){Canal 1}
\put(60,20){Canal 2}
\put(68,65){Canal 3}
\put(42,17){$P_{12}$}
\put(46,58){$P_{13}$}
\put(75,36){$P_{23}$}
\put(20,24){\line(0,1){15}}
\put(23,32){$\displaystyle s_{1}$}
\put(55,39){\line(0,1){30}}
\put(58,58){$\displaystyle s_{3}$}
\put(55,17){\line(0,1){20}}
\put(58,27){$\displaystyle s_{2}$}
\end{overpic} 
\caption{Junction - Particular case  of the straight channel.  Parameters are $\theta=\phi=0$ and $s_{1}=1$, $s_{2}=1.5$ and $s_{3}=2$.}\label{fig:straight_channel}
\end{figure}

Now, we consider the case when $\theta=\phi=0$. There is a natural way to define points $P_{12}$ and $P_{13}$,  see Fig.\ref{fig:straight_channel}. By symmetry, the $y$-coordinate of $P_{23}$ should be set to $0$ but the $x$-coordinate is undetermined.
We fix $P_{23,x}=s_1$, thus
\begin{equation*}
P_{12}\left(
\begin{array}{c}
0  
\smallskip
 \\
-s_{1}
\end{array}
\right), 
\,
P_{13}\left(
\begin{array}{c}
0
\smallskip
   \\
  s_{1}
\end{array}
\right), 
\,
P_{23}\left(
\begin{array}{c}
s_{1}
\smallskip
   \\
0
\end{array}
\right)
\end{equation*}
and  equations \eqref{eq:Conservation} reduce to:
\begin{equation*}
\left\{
\begin{aligned}
&-2  h_{1}  v_{1} + h_{2} v_{2} + h_{3} v_{3}  = 0, \\
& - 2\left(h_{1}  v_{1}^2    +\frac 1 2 g h_{1}^2\right) +
\left(h_{2}  v_{2}^2    +\frac 1 2 g h_{2}^2
\right) + \left(h_{3}  v_{3}^2     +\frac 1 2 g h_{3}^2
\right)=0, \\
&- \frac 1 2 g h_{2}^2
 + \frac 1 2 g h_{3}^2
=0. \\
\end{aligned}
\right.
\end{equation*}

Note that, with this configuration we can extend our construction to the case of a single channel with a varying cross section.

\section{Numerical scheme for shallow-water equations complemented with junction conditions}
\label{sec:NumScheme}

In this section, we couple a standard final volume scheme for the shallow water equations along each channel with the numerical flux consistent with our junction conditions \eqref{eq:cond_riemann_pb}-\eqref{eq:Conservation}. 
 
\subsection{One dimensional finite volume scheme} 
For the sake of simplicity, we will suppose that each canal has the same length, discretized with a uniform grid. Then, the computational domain in each canal is defined by the finite interval $[0,L]$, which is divided in $M$ equal cells, of amplitude $\Delta x = L/M$. The cell centers are given by $x_j=(j-\tfrac12)\Delta x$, $j=1,\ldots, M$, and the cell average of the numerical solution at time $t$ in the $j$-cell is defined as
\begin{equation}
U_j(t) = \frac{1}{\Delta x}\int_{x_{j-\frac 12}}^{x_{j+\frac 12}} U(x,t) \ dx,
\end{equation}
with appropriate boundary conditions for $U_0(\cdot)$ in channel 1 and $U_{M+1}(\cdot)$ on the two outgoing channels 2 and 3.
The system is evolved until the final time $T$, with time step $\Delta t$. We denote by  $U^n_j$ the approximate value for the average of $U$ in cell $j$ at the discrete time $t_{n}=n \Delta t$.
Hence, the finite volume approximation of system \eqref{eq:shallow_water_homog} can be written under the form
\begin{equation}\label{eq:eulero}
\frac{U^{n+1}_j-U^n_j}{\Delta t} = - \frac{1}{\Delta x}\left(\hat F^n_{j+\frac 12}-\hat F^n_{j-\frac 12}\right),
\end{equation}
where 
\begin{equation}\label{eq:numflux}
\hat F^n_{j-\frac 12} = F\left(U^n_{j-1},U^n_j\right),
\end{equation}
with $F(\cdot,\cdot)$ a proper numerical flux. We apply a Godunov type numerical flux \cite{Leveque,Toro}, computing the exact intermediate state for the Riemann problem between the two cells defined by $U_l=U^n_{j-1}$ and $U_r=U^n_j$ for $j=1,\ldots,M$. 

From now on we shall add a canal index and we use the notation $U^n_{j,k}$, $j=1,\ldots, M$,  to indicate the numerical solution along canal $k=1,2,3$ computed at time $t^n$, in $x_{j,k}$. 
The time step $\Delta t$ is fixed to satisfy the stability condition 
\begin{equation}\label{eq:cfl}
\Delta t\leq \frac{\Delta x}{\maxd_{k=1,2,3}\maxd_{1 \leq j \leq M}{\max{\{|\lambda_{1} (U_{j,k}^n)|,|\lambda_{2}(U_{j,k}^n)| \}}}},
\end{equation}
where $\lambda_{1}$ and $\lambda_{2}$ are the eigenvalues defined in \eqref{eq:eigenvalues}.

\subsection{Junction conditions and coupling with the finite volume 1D scheme}
 Let us now explain how we insert the junction conditions  \eqref{eq:cond_riemann_pb}-\eqref{eq:Conservation}  in the  finite volume numerical scheme
\eqref{eq:eulero}.  Let us write scheme \eqref{eq:eulero} in channel $k$, $k=1,2,3$ under the following form:
\begin{equation}\label{eq:eulerok}
\frac{U^{n+1}_{j,k}-U^n_{j,k}}{\Delta t} = - \frac{1}{\Delta x}\left(\hat F^n_{j+\frac 12,k}-\hat F^n_{j-\frac 12,k}\right).
\end{equation}
In a canal network, the extreme point $x_{1,k}$ or $x_{M,k}$ can be either a \textit{boundary point} or a \textit{junction point} connected with other canals. In our setting $x_{M,1}$, $x_{1,2}$ and $x_{1,3}$ are junction points, while $x_{1,1}$, $x_{M,2}$ and $x_{M,3}$ are boundary points. At the boundary points of the network, the numerical tests use homogeneous Neumann conditions, but other boundary conditions can naturally be used.

Let $U_{k}=(h_k,v_k)$, $k=1,2, 3$ be the solution to system \eqref{eq:cond_riemann_pb}-\eqref{eq:Conservation} with $U^*_{1}=U^n_{M,1}$, $U^*_2=U^n_{1,2}$, $U^*_3=U^n_{1,3}$ where
$U^n_{M,1}$, $U^n_{1,2}$ and $U^n_{1,3}$ are the values computed by the 1D scheme along the channels. Then, at the junction points we impose $\hat F^n_{M+1/2,1}=f(U_1)$, $\hat F^n_{1/2,2}=f(U_2)$ and $\hat F^n_{1/2,3}=f(U_3)$.

\subsection{Solving the non-linear system at the junction}\label{sec:Sol_System}
Now, let us study the solutions to the nonlinear system at the junction. 
We recall the notations used in Sec.\ref{sec:Junction} and denote by
$h_{k}^*$ and  $v_{k}^*$ the approximate values of $h$ and $v$ near the junction at channel $k=1,2,3$ given by the 1D numerical scheme, see also Sec.\ref{sec:NumScheme} for their exact definition.
Let $\Omega$ be the open set of  admissible states,
 $\Omega=\{h_k\in\mathbb{R}^+_{*}, v_k\in\mathbb{R} , |v_{k}| < \sqrt{gh_{k}}, k=1,2,3\}$.
The approximate values $h_{k}$ and  $v_{k}$ are then obtained solving the non-linear system \eqref{eq:Conservation} with $v_k$ given by  \eqref{eq:cond_riemann_pb}.
Let us denote
$$ X=\left(\begin{array}{c} h_{1} \\ h_{2} \\ h_{3} \\ v_{1}\\ v_{2}\\ v_{3} \end{array}\right), \; X^*= \left(\begin{array}{c} h_{1}^* \\ h_{2}^* \\ h_{3}^* \\ v_{1}^*\\ v_{2}^*\\ v_{3}^* \end{array}\right) 
= \left(\begin{array}{c} h_{1,M} \\ h_{2,1} \\ h_{3,1} \\ v_{1,M}\\ v_{2,1}\\ v_{3,1} \end{array}\right).$$
We can rewrite the system under the following form:
\begin{equation}\label{system}
  \Psi
  \left( 
  X;X^*
  \right)=0
  \end{equation}
where $\displaystyle \Psi :\Omega  \times  \Omega \to \mathbb{R}^6$.
In general,  existence and uniqueness results for solutions to non linear systems are difficult to prove.
 
Assume that we have solved the system up to $t=t^n$, this gives the solution $X^{n,*}$ which faces the junction.
Suppose that we have found a solution $X^n$ such that $\Psi(X^n;X^{n,*})=0$.
If we can prove that $\textrm{ Det } D\Psi \left(X^n;X^{n,*}\right)\neq 0$, where $D\Psi$ denotes the Jacobian with respect to  the first argument, then there exists a unique $X=X(X^*)$, for $\|X^*-X^{n,*}\|<\epsilon$, with $\epsilon$ small enough, such that $\Psi(X;X^{*})=0$.
Therefore, if the flow is smooth, one can find $\Delta t$ small enough such that  $\|X^{n+1,*}-X^{n,*}\|<\epsilon$ and the implicit function theorem guarantees that there exists a unique solution $X^{n+1}$ such that $\Psi(X^{n+1};X^{n+1,*})=0$. So, the procedure can be iterated provided 
one can prove at each step that $\textrm{ Det } D\Psi \left(X^n;X^{n,*}\right)\neq 0$.

In the particular case when all waves are rarefactions, the 
relations \eqref{eq:cond_riemann_pb} become
\begin{equation}\label{junction_rar}
\left\{
\begin{aligned}
v_{1}+ 2 \sqrt{gh_{1}} =v_{1}^*+2 \sqrt{gh_{1}^*},\\
v_{2}- 2 \sqrt{gh_{2}} =v_{2}^*-2 \sqrt{gh_{2}^*},\\
v_{3}- 2 \sqrt{gh_{3}} =v_{3}^*-2 \sqrt{gh_{3}^*},\\
\end{aligned}
\right.
\end{equation}
and it is clear that the Jacobian $D\Psi$ does not depend on the data $X^*$. 
Thus, starting from a set of data $X^*$ and a solution $X$ such that $\Psi(X;X^*)=0$, 
once one can prove that 
\begin{equation}\label{condition}
  \textrm{ Det } D \Psi \left(X;X^{*}\right) \neq 0, 
  \text{ for all }X^{*},X \in \Omega  \times  \Omega, 
  \end{equation}
  the solution exists at each time step.
Note that, for the steady solution $h=const.$ and $v=0$ one has
$X^*=(h,h,h,0,0,0)^T$ and $\Psi(X;X^*)=0$ for $X=X^*$, thus there exists at least one case for which $\Psi(X;X^*)=0$.

\begin{example} In the case when the canals are orthogonal to the sides of the triangle, we have
\begin{equation*}\label{function}
  \Psi
  \left( 
  X;X^*
  \right)=\left(\begin{array}{c} 
  s_{1}h_{1}v_{1}-s_{2} h_{2}v_{2}-s_{3} h_{3}v_{3} \\ 
  \medskip 
  h_{1}v_{1}^2+\frac 1 2 gh_{1}^2- h_{2}v_{2}^2-\frac 1 2 g h_{2}^2\\ 
    \medskip 
  h_{1}v_{1}^2+\frac 1 2 gh_{1}^2- h_{3}v_{3}^2-\frac 1 2 g h_{3}^2 \\
    \medskip 
  v_{1}+ 2 \sqrt{gh_{1}} -v_{1}^*-2 \sqrt{gh_{1}^*}\\
      \medskip 
    v_{2}- 2 \sqrt{gh_{2}} -v_{2}^*+2 \sqrt{gh_{2}^*}\\ 
   \medskip 
 v_{3}- 2 \sqrt{gh_{3}} -v_{3}^*+2 \sqrt{gh_{3}^*}
  \end{array}\right).
\end{equation*}
Tedious algebra gives,
\begin{equation*}
   \begin{aligned}
  \textrm{ Det } D \Psi \left(X;X^*\right) 
  &=      \textrm{ Det }
   \left( \begin{array}{ccc} 
s_{1} (v_{1}-\sqrt{gh_{1}})&-s_{2} (v_{2}+\sqrt{gh_{2}})&-s_{3} (v_{3}+\sqrt{gh_{3}})\\
(v_{1}-\sqrt{gh_{1}})^2&- (v_{2}+\sqrt{gh_{2}})^2&0\\
(v_{1}-\sqrt{gh_{1}})^2&0&- (v_{3}+\sqrt{gh_{3}})^2\\
  \end{array}\right)\\
&=\lambda_1 \lambda_2 \lambda_3\left( s_{1}\lambda_2 \lambda_3-s_{2}\lambda_1 \lambda_3-s_{3}\lambda_1\lambda_2\right).
   \end{aligned}
  \end{equation*}
We can therefore conclude that since we are in the sub-critical case, for which $$
\lambda_1 = v_{1}-\sqrt{g h_{1}}<0, \quad \lambda_2 = v_{2} + \sqrt{g h_{2}}>0, \quad \lambda_3= v_{3} + \sqrt{g h_{3}}>0,
$$
we have 
  $$\textrm{ Det }  D \Psi \left(X;X^{*}\right) <0,   \text{ for all }X^{*},X \in \Omega, $$
  which implies condition \eqref{condition}. Thus, starting from a point for which $\Psi(X;X^*)=0$ we can prolong the solution for all of times. This coincides with the case in which one assumes the continuity of the energy.
\end{example}

\begin{example}

\textit{Case with vanishing velocities.}

Now, consider the case given by system\eqref{eq:Conservation}-\eqref{eq:cond_riemann_pb}, with only rarefaction waves. In order to simplify the expressions arising in the  computations, we introduce the following notations
 \begin{equation*}
\left\{
\begin{aligned}
\alpha_{1}&= \ell_{1} \left(\begin{array}{c} 1 \\ 0\end{array}\right)\cdot \mathbf{n}_{1}= -2s_{1},  \\
\alpha_{2}&= \ell_{2} \left(\begin{array}{c} 1 \\ 0\end{array}\right)\cdot \mathbf{n}_{2}= s_{1}+\displaystyle\frac{s_{2}\sin \theta +s_{3} \sin\phi}{\sin (\theta-\phi)} ,
  \\
\alpha_{3}&= \ell_{3} \left(\begin{array}{c} 1 \\ 0\end{array}\right)\cdot \mathbf{n}_{3}=    s_{1}-\displaystyle\frac{s_{2}\sin \theta +s_{3} \sin\phi}{\sin (\theta-\phi)} , 
 \\
\beta_{1}&= \ell_{1} \left(\begin{array}{c} 0 \\ 1\end{array}\right)\cdot \mathbf{n}_{1}=\displaystyle\frac{s_{1} \sin (\theta+\phi)-s_{2}\sin \theta -s_{3} \sin\phi}{\sin \phi\sin \theta} , 
 \\
\beta_{2}&= \ell_{2} \left(\begin{array}{c} 0 \\ 1\end{array}\right)\cdot \mathbf{n}_{2}=\displaystyle-\frac{s_{2}\cos \theta +s_{3} \cos\phi}{\sin (\theta-\phi)} -
\frac{s_{1} \cos \phi-s_{2}}{\sin \phi} ,  \\
\beta_{3}&= \ell_{3} \left(\begin{array}{c} 0 \\ 1\end{array}\right)\cdot \mathbf{n}_{3}=\displaystyle \frac{s_{2}\cos \theta +s_{3} \cos\phi}{\sin (\theta-\phi)} -\frac{s_{1}\cos \theta-s_{3}}{\sin \theta},  \\
\gamma_{2}&= \ell_{2} \left(\begin{array}{c} \cos\phi \\ \sin\phi\end{array}\right) \cdot \mathbf{n}_{2}=\alpha_{2} \cos\phi+\beta_{2} \sin \phi,\\
\gamma_{3}&=\ell_3 \left(\begin{array}{c} \cos\theta \\ \sin\theta\end{array}\right)\cdot \mathbf{n}_{3} =\alpha_{3} \cos\theta+\beta_{3} \sin \theta,
\end{aligned}
\right.
\end{equation*}
such that 
 \begin{equation*}
 \Psi
  \left( 
  X;X^*
  \right)=\left(
\begin{aligned}
& \alpha_1 h_1 v_{1}   + \gamma_2 h_2 v_{2}   + \gamma_3 h_3 v_{3}  
\medskip\\
&  
 \alpha_1\left( h_{1}v_{1}^2    +\frac 1 2 g h_{1}^2 \right) +
\left(\gamma_2 h_{2}v_{2}^2  \cos\phi    +\frac{\alpha_2}{2}   g h_{2}^2
\right) 
+\left(\gamma_3 h_{3}v_{3}^2 \cos\theta    +\frac{\alpha_3}{2} g h_{3}^2
\right) \\
& \frac 1 2  \beta_{1}g h_{1}^2
 +
\left(\gamma_2 h_{2}v_{2}^2  \sin\phi     +\frac{\beta_2}{2}g h_{2}^2
\right)
 +\left(\gamma_3 h_{3}v_{3}^2  \sin\theta    +\frac{\beta_3}{2} g h_{3}^2
\right) \\
&v_{1}+ 2 \sqrt{gh_{1}} -v_{1}^*-2 \sqrt{gh_{1}^*}\\
&v_{2}- 2 \sqrt{gh_{2}} -v_{2}^*+2 \sqrt{gh_{2}^*}\\
&v_{3}- 2 \sqrt{gh_{3}} -v_{3}^*+2 \sqrt{gh_{3}^*}\\
\end{aligned}
\right).
\end{equation*}

Since explicit computations are too difficult, we restrict ourselves to the particular case when solutions with  vanishing velocities $v_{1}^n=v_{2}^n=v_{3}^n=0$ at the junction exist, that is to say 
$$ X^n=\left(\begin{array}{cccccc} h_{1}^n & h_{2}^n & h_{3}^n & 0 & 0 & 0 \end{array}\right).$$ 

We already note that $\ds \mathbf{n}_{1}+\mathbf{n}_{2}+\mathbf{n}_{3} =0$, if $\ds (h_{1}^*, h_{2}^*, h_{3}^*, v_{1}^*, v_{2}^*, v_{3}^*)=(h,h,h,0,0,0)$, $X^n=\ds (h_{1}, h_{2}, h_{3}, v_{1}, v_{2}, v_{3})=(h,h,h,0,0,0)$ is a   trivial solution   of system \eqref{eq:Conservation}-\eqref{eq:cond_riemann_pb}.

We prove that in this case $\textrm{ Det } D\Psi \left(X;X^{*} \right)\neq 0$. In fact,

  \begin{equation*}
D\Psi \left(X;X^{*}\right) =\left( \begin{array}{ccc|ccc} 
 0&0&0&\alpha_{1} h_{1}&\gamma_{2} h_{2}& \gamma_{3} h_{3}\\
 \alpha_{1} g h_{1}& \alpha_{2} g h_{2}& \alpha_{3} g h_{3}&0&0&0\\
\beta_{1} g h_{1} &\beta_{2} g h_{2}&\beta_{3} g h_{3}&0&0&0\\
\hline
\frac{\sqrt{g}}{\sqrt{h_{1}}}&0&0&1&0&0\\
0&-\frac{\sqrt{g}}{\sqrt{h_{2}}}&0&0&1&0\\
0&0&-\frac{\sqrt{g}}{\sqrt{h_{3}}}&0&0&1
  \end{array}\right)
  \end{equation*}
  and 
   \begin{equation*}
   \begin{aligned}
   \textrm{ Det } D \Psi \left(X;X^{*}\right) &= g^{5/2}
      \textrm{ Det }
   \left( \begin{array}{ccc} 
-\alpha_{1} \sqrt{ h_{1}}&\gamma_{2} \sqrt{ h_{2}}& \gamma_{3} \sqrt{ h_{3}}\\
 \alpha_{1}  h_{1}& \alpha_{2}  h_{2}& \alpha_{3}  h_{3}\\
\beta_{1}  h_{1} &\beta_{2}  h_{2}&\beta_{3} h_{3}
  \end{array}\right)\\
    &= g^{5/2}\sqrt{h_{1}h_{2} h_{3}}\bigl(-\alpha_{1}(\alpha_{2}\beta_{3}-\alpha_{3}\beta_{2})\sqrt{h_{2} h_{3}}-\gamma_{2}(\alpha_{1}\beta_{3}-\alpha_{3}\beta_{1})\sqrt{h_{1} h_{3}} \\ 
    & \qquad \qquad\qquad  +\gamma_{3}(\alpha_{1}\beta_{2}-\alpha_{2}\beta_{1})\sqrt{h_{1} h_{2}} \bigr)\\
  &=2 g^{5/2}\sqrt{h_{1}h_{2} h_{3}} (\sqrt{h_{1}h_{2}}s_{3}+ \sqrt{h_{1}h_{3}}s_{2} + \sqrt{h_{2}h_{3}}s_{1})\times \Bigl( \frac{4 s_{2}s_{3} }{\sin( \theta-\phi)}+\\
  & \qquad \qquad\qquad \frac{(s_{1}\sin( \theta-\phi) +s_{3}\sin(\phi) - s_{2}\sin(\theta))^2}{\sin( \theta-\phi)\sin(\phi)\sin(\theta)}\Bigr)
   \end{aligned}
  \end{equation*}
  
  We notice that this expression is the product of two terms: the first one, $\displaystyle 2 g^{5/2}\sqrt{h_{1}h_{2} h_{3}} (\sqrt{h_{1}h_{2}}s_{3}+ \sqrt{h_{1}h_{3}}s_{2} + \sqrt{h_{2}h_{3}}s_{1})$ is always positive since $h_{1}>0, \, h_{2}>0, \,h_{3}>0$. 
  
  The second term  depends only on $s_{1}, s_{2}, s_{3}, \theta, \phi$, that is to say on the triangle geometry. This second term vanishes iff the triangle is degenerate, see Remark~\ref{alignment}.
  
  Therefore, excluding the case of a degenerate triangle,   
  $$ \textrm{ Det }  D \Psi \left(X;X^{*}\right) \neq 0.$$
  
  \end{example}
  
\section{Numerical tests}\label{Sec:NumTests}
In order to evaluate the effectiveness of our method we perform various tests consisting of
a subcritical wave propagating across junctions of different geometries. For that purpose, we will use the numerical scheme presented in the previous section.
First of all we check the numerical convergence of the scheme on the whole network under grid refinement. Then, we compare numerical simulations on a network with $\theta=\phi=0$,
with simulations on a single canal, to show the consistency of the junction conditions with the traditional one-dimensional Riemann solver. Subsequently, we increase $\theta$ and $\phi$ in order to enhance the influence of the angles in our junction conditions. We also investigate numerically the case when water flows out canals 1 and 3  and pours into canal 2. Finally, we compare the dynamics of the 1D solver with the numerical solution obtained with a two-dimensional code.

\subsection{Convergence of the 1D numerical scheme}\label{test:Conv1D}
\begin{figure}[htbp!]
\includegraphics[scale=0.2]{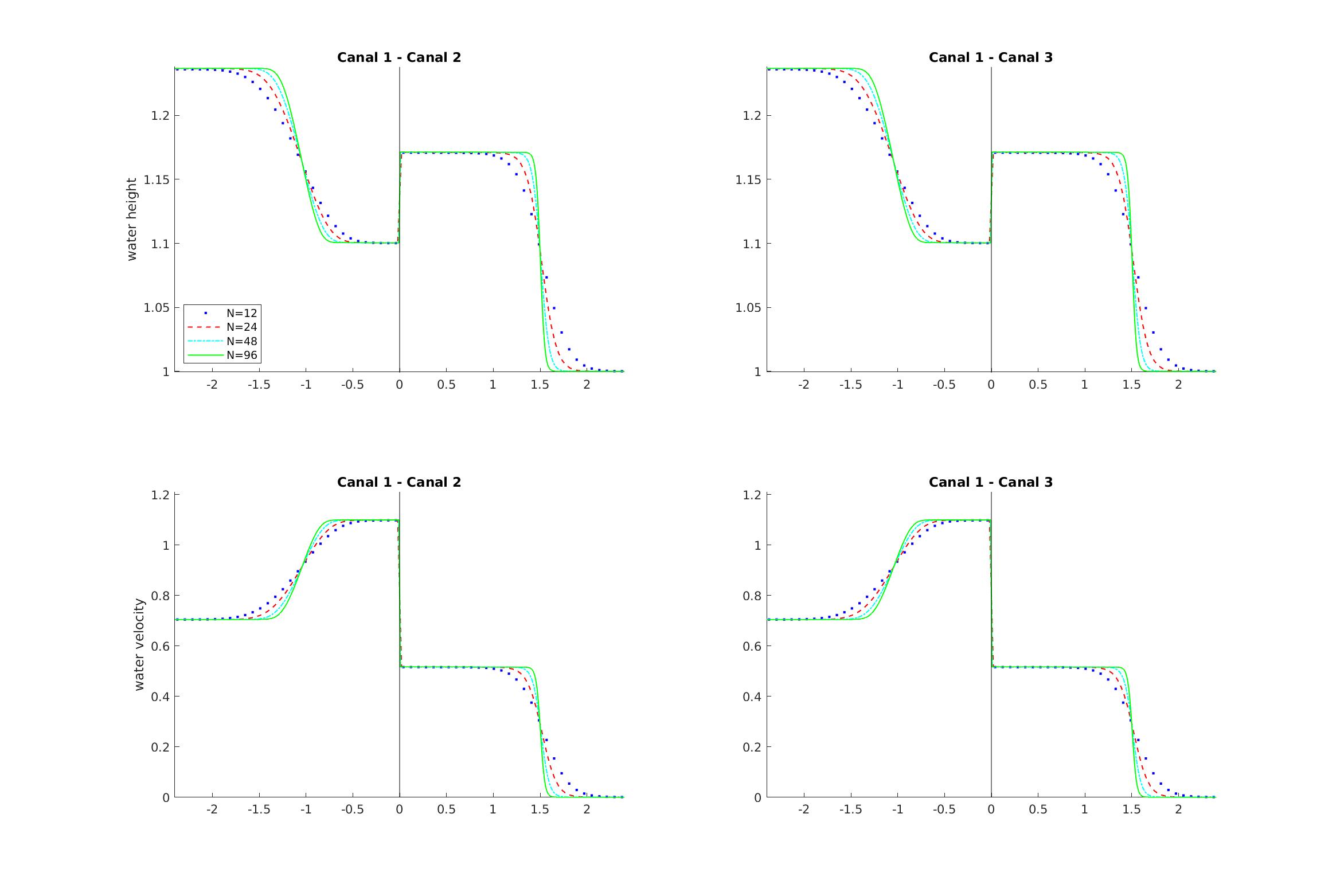}
\caption{Test \ref{test:Conv1D}: numerical convergence at time $T=0.9$, after the shock has reached the junction. $N=12,24,48,96$. On top: water height as a function of $x$; on bottom: water velocity as a function of $x$. On the left: canal 1 and 2; on the right: canal 1 and 3. Initial data are given in \eqref{eq:initdata1}.}\label{fig:Conv1Ddopo}
\end{figure}

In this test case we check the numerical convergence of the 1D scheme \eqref{eq:eulero} coupled with system \eqref{eq:Conservation}-\eqref{eq:cond_riemann_pb} at the junction under grid refinement. We set the geometry parameters 
$s_1=s_2=s_3=1$,  $\theta=-\phi=\pi/6$, $L_{1}=L_{2}=L_{3}=3.5$, where $L_k$ is the length of channel $k$. 
The initial data in the three channels is
\begin{equation}\label{eq:initdata1}
\begin{array}{l}
v_k(x,t=0)=0, \quad k=1,2,3,  \, x \in [0,L_{k}], 
\\
h_2(x,t=0)=h_3(x,t=0)=1, 
h_1(x,t=0) = \left\{\begin{array}{ll}
1.5 & x\leq L_1/2
\\
1  & x> L_1/2.
\end{array}\right.
\end{array}
\end{equation}
We expect the formation of a rarefaction wave propagating backwards on channel 1 and a shock crossing the junction and travelling with positive speed along the two outgoing canals. In Figure \ref{fig:Conv1Ddopo} we show the solution obtained after the water wave has reached the junction. The dynamic involves the three canals and we observe the convergence of the numerical solution under grid refinement. The number of grid points on each channel is $N=12,24,48,96$. We observe that the solution of system  \eqref{eq:Conservation}-\eqref{eq:cond_riemann_pb} 
at the junction does not depend on the grid parameter $N$, proving the consistency of the Riemann solver. We observe the formation of a stationary shock at the junction.

\subsection{Comparison with the solution on a single canal}\label{test:SingleCanal}
Here we study the influence of  our junction conditions involving the angles $\theta$ and $\phi$ on the solution.

We set $s_1=1$ and $s_2=s_3=1/2$, $L_{1}=L_{2}=L_{3}=5$ and the initial data are the same as in \eqref{eq:initdata1}.
In all tests we fix the grid parameters $\Delta x=0.01$ on each channel and $\Delta t$ to satisfy \eqref{eq:cfl}.

\begin{figure}[h!]
\begin{overpic}
 [width=1.0\textwidth]{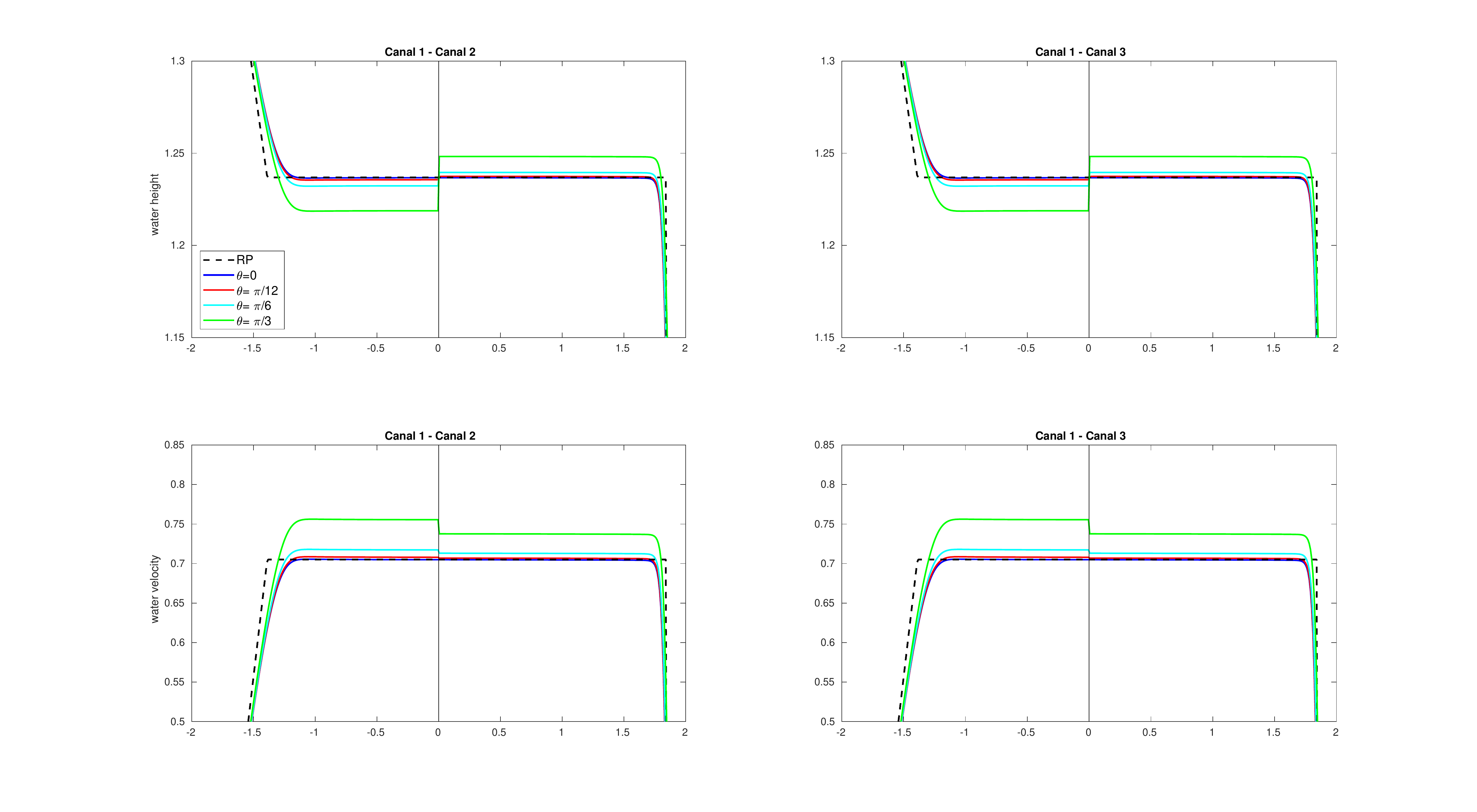}
  \end{overpic}
\caption{Test \ref{test:SingleCanal}: comparison of various symmetric geometries with the solution on a single canal at $T=0.5$: $s_1=1$ and $s_2=s_3=1/2$;  $\phi = -\theta$ and $\theta=0,\pi/12,\pi/6,\pi/3$. Black dashed line: 1D exact shallow-water on a single channel; Blue solid line: numerical solution for $\theta=\phi=0$. On the left: zoom on the transition between canal 1 and 2; on the right: zoom on the transition between canal 1 and 3. Initial data are given in \eqref{eq:initdata1}.}\label{fig:SingleCanalZoom_sim}
\end{figure}

We first consider $\theta=\phi=0$ which corresponds to a single channel. 
As expected for $\theta=\phi=0$ and $s_1=s_2+s_3$ the solution of our algorithm coincides with the solution computed on a single canal, i.e. without the junction, see blue line and black dashed line on Figure \ref{fig:SingleCanalZoom_sim}.
\\
Next, 
we change $\theta$ and $\phi$ to study the influence of the angles on the dynamics.
Specifically, in Figure \ref{fig:SingleCanalZoom_sim} the angles vary symmetrically with $\theta=-\phi$ and $\theta=0,\pi/12,\pi/6,\pi/3$. We observe that the symmetry of the configuration is preserved and that the solution varies monotonically increasing the angles and  moving away from the single channel profile. Note that to see the dependence of the solution on the angles, it is essential to include an angle dependence in the junction condition.

\begin{figure}[h!]
\begin{overpic}
 [width=1.0\textwidth]{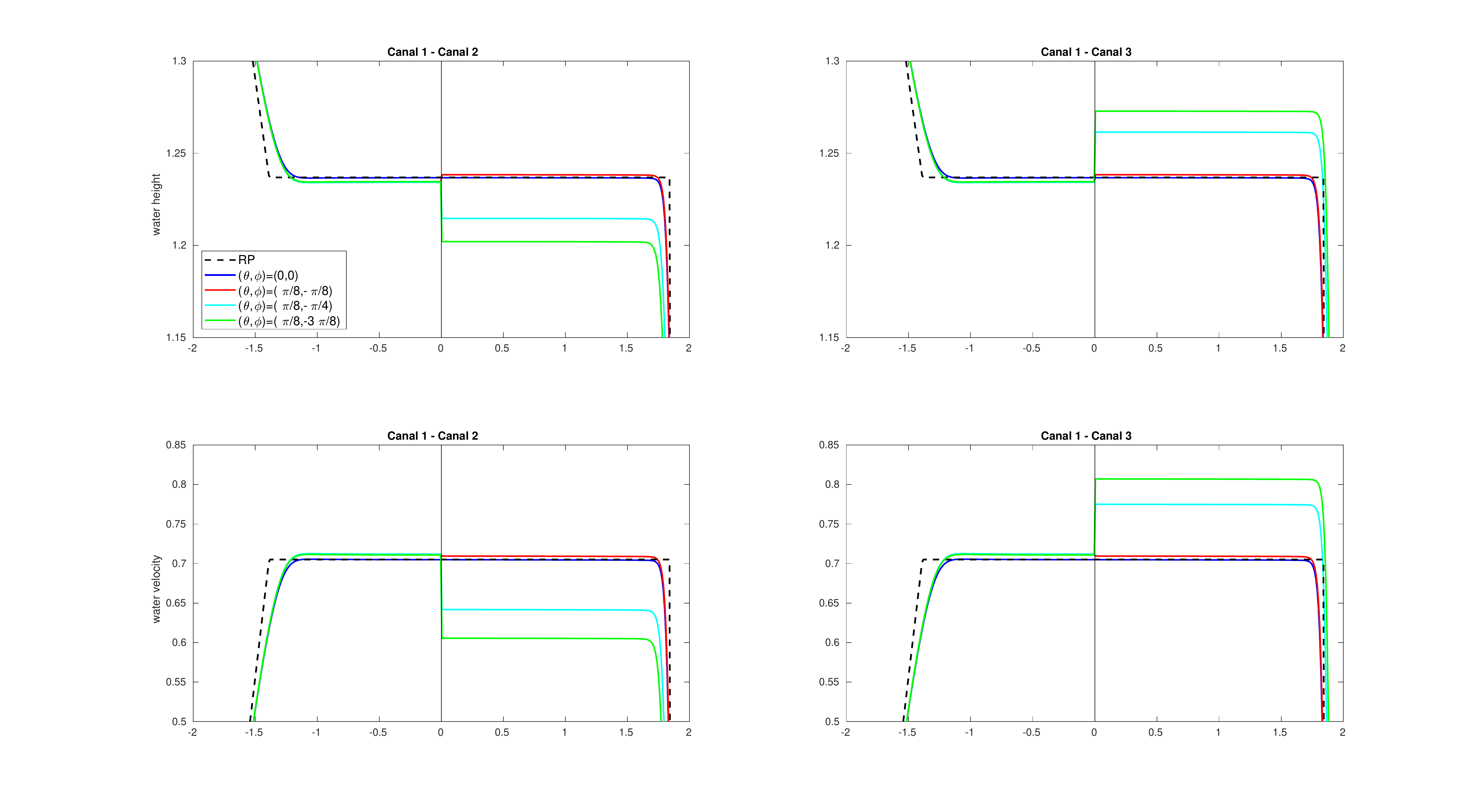}
  \end{overpic}
\caption{Test \ref{test:SingleCanal}: comparison of various triangle asymmetric geometries with different angles with the solution on a single canal at $T=0.5$. $s_1=1$ and $s_2=s_3=1/2$, $\theta=\pi/8$ fixed and $\phi=-\pi/8,-\pi/4,-3\pi/8$. Black dashed line: 1D exact shallow-water on a single channel; Blue solid line: numerical solution for $\theta=\phi=0$. On the left: zoom on the transition between canal 1 and 2; on the right: zoom on the transition between canal 1 and 3. Initial data are given in \eqref{eq:initdata1}.
}\label{fig:SingleCanal_asim1}
\end{figure}

In Figure \ref{fig:SingleCanal_asim1} we study the influence of a non symmetric variation of the angles: we fix $\theta=\pi/8$ and vary $\phi=-k\pi/8$ with $k=0,1,2,3$. We observe that the symmetry of the solutions of the two outgoing channels is lost and that  water meets more resistance as channel 2 becomes more and more bent. So, the water level decreases in channel 2 and increases in channel 3. The solution on channel 1 does not change because the total lumen of the outgoing channels remains the same.

\begin{figure}[htbp!]
\begin{overpic}
 [width=1.0\textwidth]{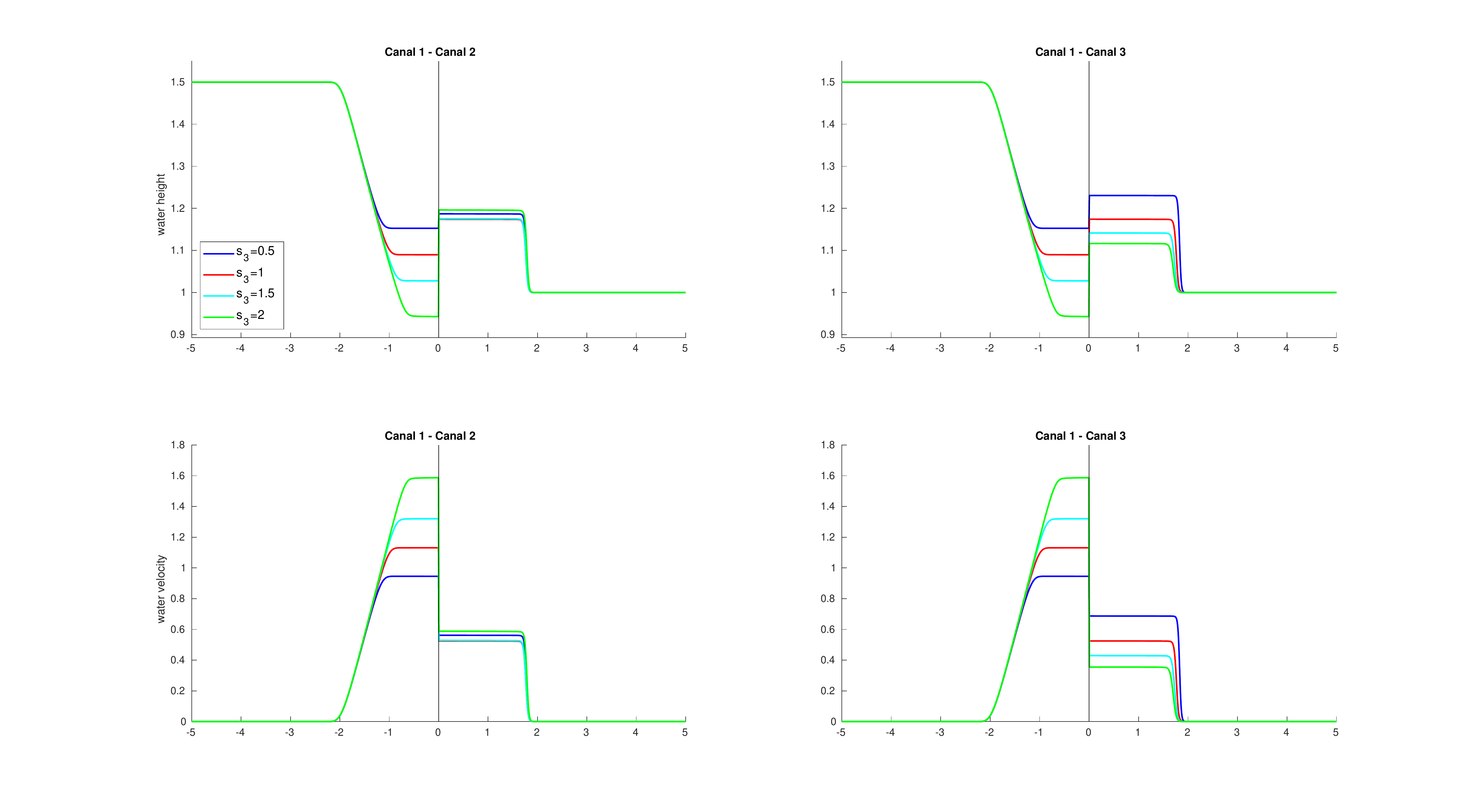}
  \end{overpic}
\caption{Test \ref{test:SingleCanal}: comparison of various geometries with different sections at $T=0.5$.  $\theta=-\phi=\pi/4$, $s_1=s_2=1$ and $s_3=0.5,1,1.5,2$. 
On the left: zoom on the  transition between canal 1 and 2; on the right: zoom on the  transition between canal 1 and 3. Initial data are given at Eq.\eqref{eq:initdata1}.
}
\label{fig:test_sizes}
\end{figure}

Finally, we compare the 1D solution fixing the two angles and varying the channel sizes. Specifically, we fix $\theta=-\phi=\pi/4$, $s_1=s_2=1$ and consider $s_3=0.5,1,1.5, 2$.
In Figure \ref{fig:test_sizes} we observe that, as the section of channel 3 increases,  the water height decreases in canals 1 and 3 and increases in canal 2. The water velocity in channels 2 and 3 follows the same behavior, while the velocity in channel 1 increases. However, the dynamics in channel 2 does not vary significantly.

\subsection{Merging canals}
\label{sec:test_2to1}
\begin{figure}[htbp!]
\begin{overpic}
 [width=1.0\textwidth]{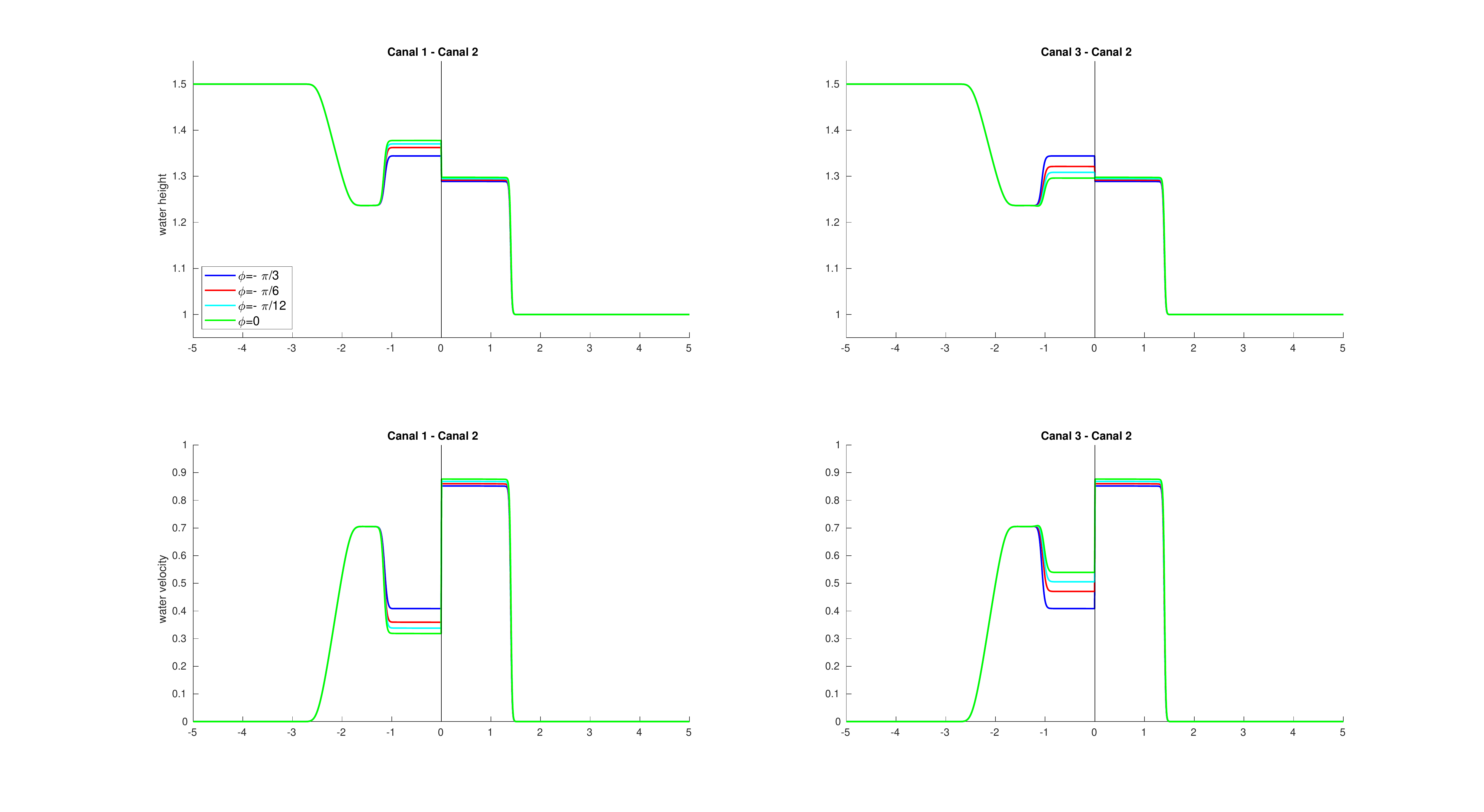}
  \end{overpic}
\caption{Test \ref{sec:test_2to1}: comparison of the solution for various geometries in the case of a merging at $T=0.5$, $s_1=s_2=s_{3}=1$,  $\theta=\pi/3$ and $\phi=-\pi/3,-\pi/6,-\pi/12,0$.
On the left: zoom on the  transition between canal 1 and 2; on the right:   zoom on the  transition between canal 3 and 2. Initial data are given at Eq.\eqref{eq:initdata2}.
}\label{fig:confronto1d_2in1}
\end{figure}

In this section we display numerical results for a 
2-to-1 or \textit{merging} junction for which the water flows from channels 1 and 3 towards channel 2. 
We set as initial data
\begin{equation}\label{eq:initdata2}
\begin{array}{l}
v_k(x,t=0)=0,  k=1,2,3,    \, x \in [0,L_{k}], 
\\
h_2(x,t=0)= 1,\  
h_i(x,t=0) = \left\{\begin{array}{ll}
1.5 & x\leq L_{i}/2
\\
1  & x> L_{i}/2.
\end{array}\right., i=1,3.
\end{array}
\end{equation}
Then, we fix $s_{1}=s_{2} =s_{3}=1$,  $\theta=\pi/3$. In Figure \ref{fig:confronto1d_2in1}, we compare the solutions obtained for $\phi=-\pi/3,-\pi/6,-\pi/12,0$. As the angle $\phi$ widens, the water height increases in channel 1  and decreases in channel 3. This asymmetry explains why the dynamic in channel 2 is almost unaffected by the angle variation, the amount of water entering remains almost constant.

\subsection{Comparison of 1D and 2D solutions}\label{sec:test_1D2D}
We compare our 1D solver with the 2D shallow water solution  \eqref{eq:saint_venant2D}.

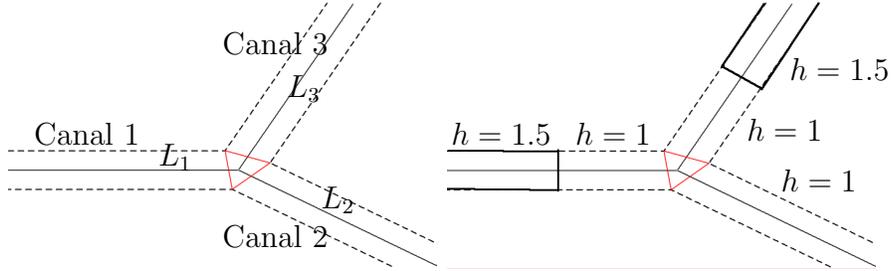
\begin{figure}[htbp!]
\begin{overpic}
[width=0.45\textwidth]{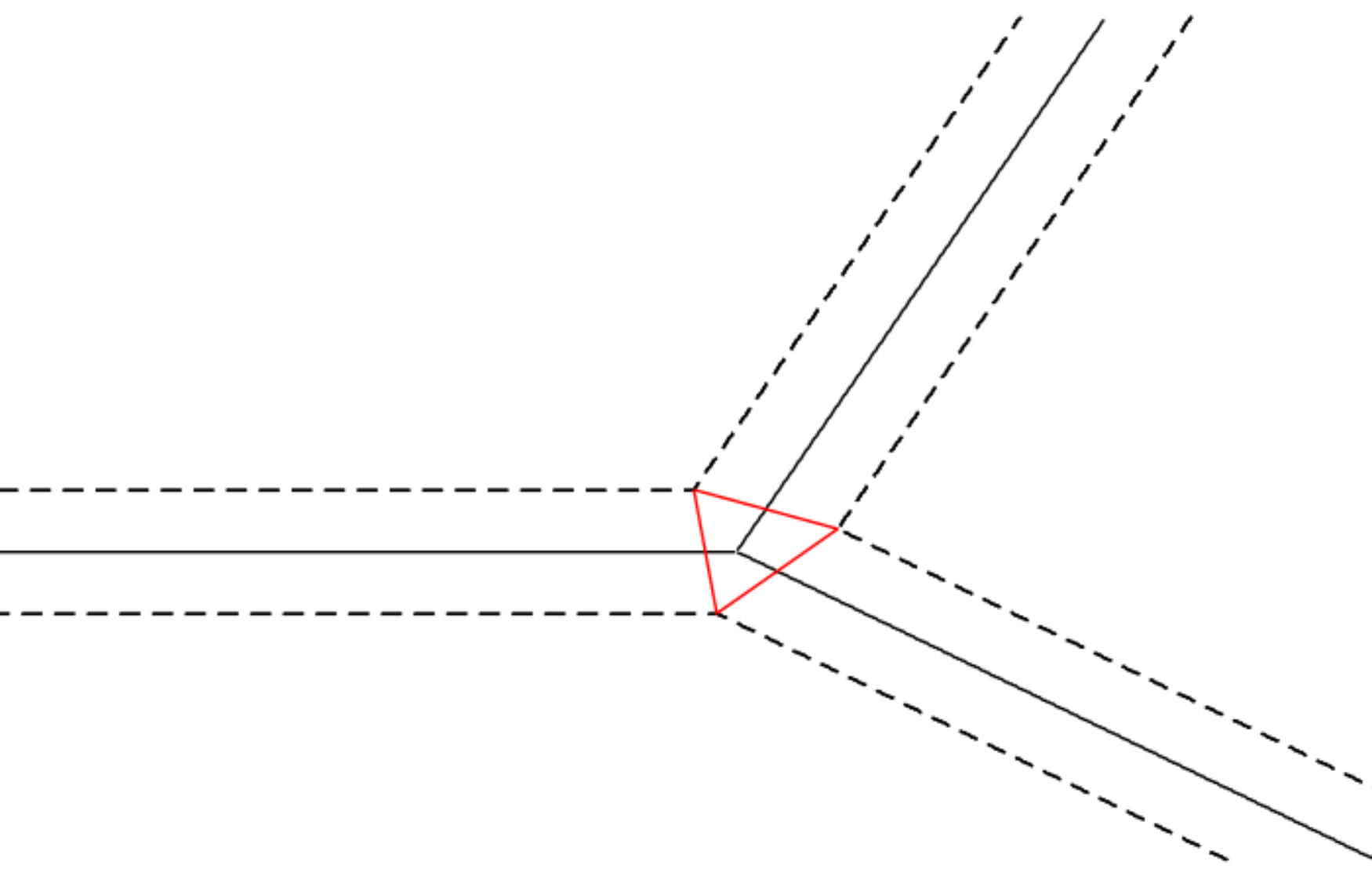}
\put(6,29){Canal 1}\put(35,24){$L_1$}
\put(50,50){Canal 3}\put(65,40){$L_3$}
\put(50,5){Canal 2}\put(73,14){$L_2$}
\end{overpic} 
\begin{overpic}
[width=0.45\textwidth]{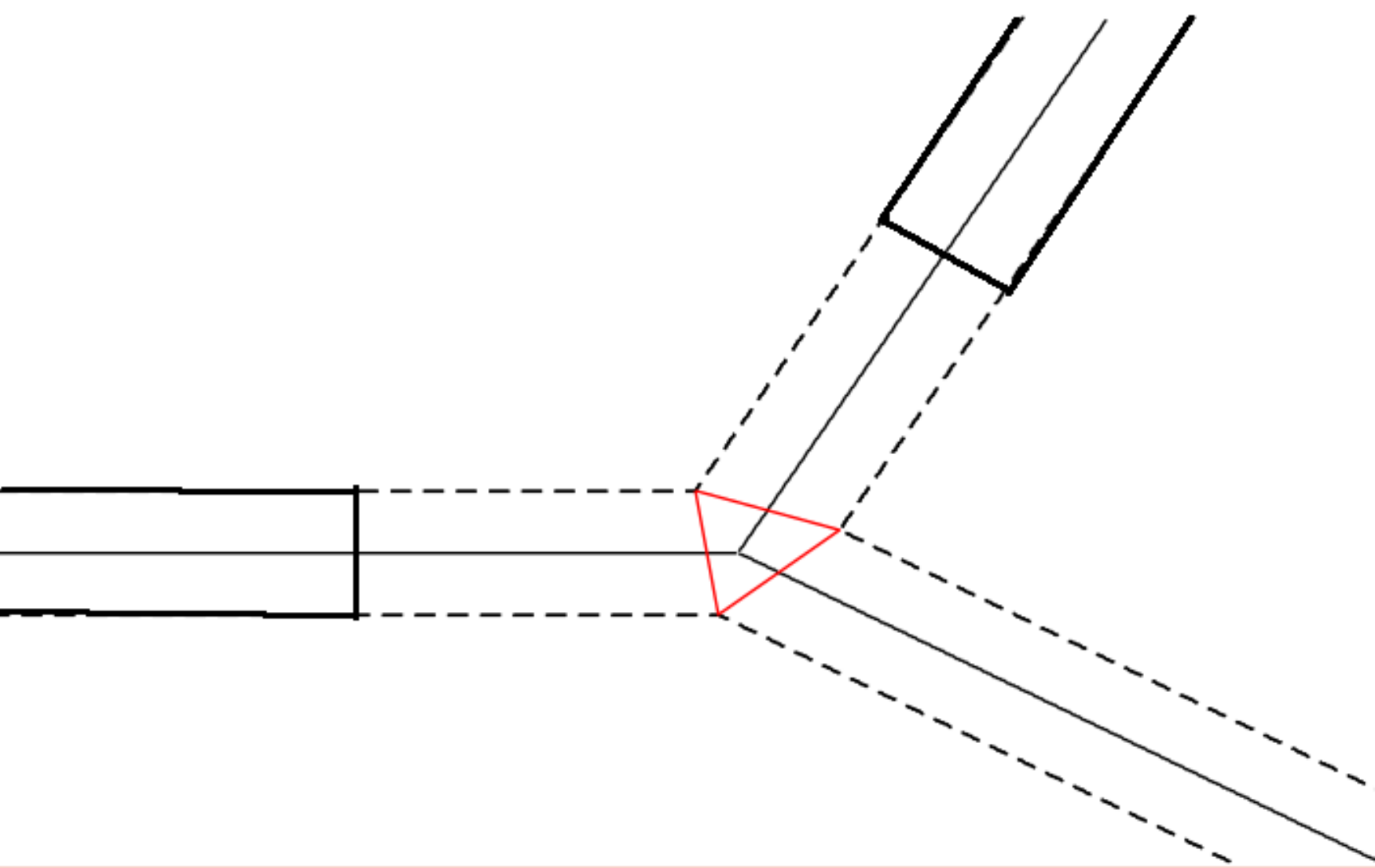}
\put(1,29){$h=1.5$}\put(30,29){$h=1$}
\put(80,44){$h=1.5$}\put(70,30){$h=1$}
\put(78,18){$h=1$}
\end{overpic} 
\caption{On the left the 2D numerical domain for the comparison with the 1D geometry. On the right the 2D initial water height in a merging junction.}\label{Fig:geom2D}
\end{figure}

The numerical solution of \eqref{eq:saint_venant2D} has been computed by the free and open source ToolBox \textit{FullSWOF2D} (Full Shallow-Water equations for Overland Flow in 2D), which is a C++ code for simulations in two dimensions \cite{FullSWOFpaper}. We compute the solution on the rectangle in Figure \ref{Fig:geom2D} with $\displaystyle (v_{x}, v_{y})\cdot \mathbf{n}=0$ on the boundary. To obtain the 2D geometry of the junction we are interested in, we use a bottom topography which is $z_{in}=-1$ within the dashed region $\Omega$ and $z_{out}=0$ in the complement $\Omega^C$. We choose the initial water height so that $h+z_{in}$ is less than $z_{out}$. In this way, the water flows only inside the dashed region while $\Omega^C$ is seen as a dry state region.

To compare the results of the 2D solution with the 1D code, the values of $h$, $hv_x$, $hv_y$ of the 2D solution are sampled on the straight lines at the middle of the channels (black solid  lines in Figure \ref{Fig:geom2D}) of lengths $L_k$, $k=1,2,3$. We compare the 1D velocity with the 2D velocity norm  $\sqrt{v_x^2+v_y^2}$.
The mesh of the two-dimensional domain contains about $26\times 10^4$ grid nodes.

We fix the canal lengths as $L_1=L_2=L_3=5$ and the simulation final time to $T=1.5$.  We set $s=s_1=s_2=s_3$, and  we compare the 1D solution obtained with junction conditions \eqref{eq:Conservation}-\eqref{eq:cond_riemann_pb}  to the 2D solution with $s=1,\ 0.5, \ 0.25.$  As $s$ decreases, the 2D configuration becomes closer to the 1D network.
Recall that, in the 1D case, if $s=s_1=s_2=s_3$, the solution of the shallow water equations complemented with junction conditions \eqref{eq:Conservation}-\eqref{eq:cond_riemann_pb}  does not depend on the value of $s$.

We consider two different cases, 
a \textit{diverging junction}, 
that it to say the case when  water flows  from the single channel 1 towards the two channels 2 and 3 and a \textit{merging junction}, when the water flows from the two channels 1 and 3 towards the single  channel 2.

\medskip 

\begin{figure}[htbp!]
 \subfloat[water height]
 {
\begin{overpic}
 [width=1.0\textwidth]{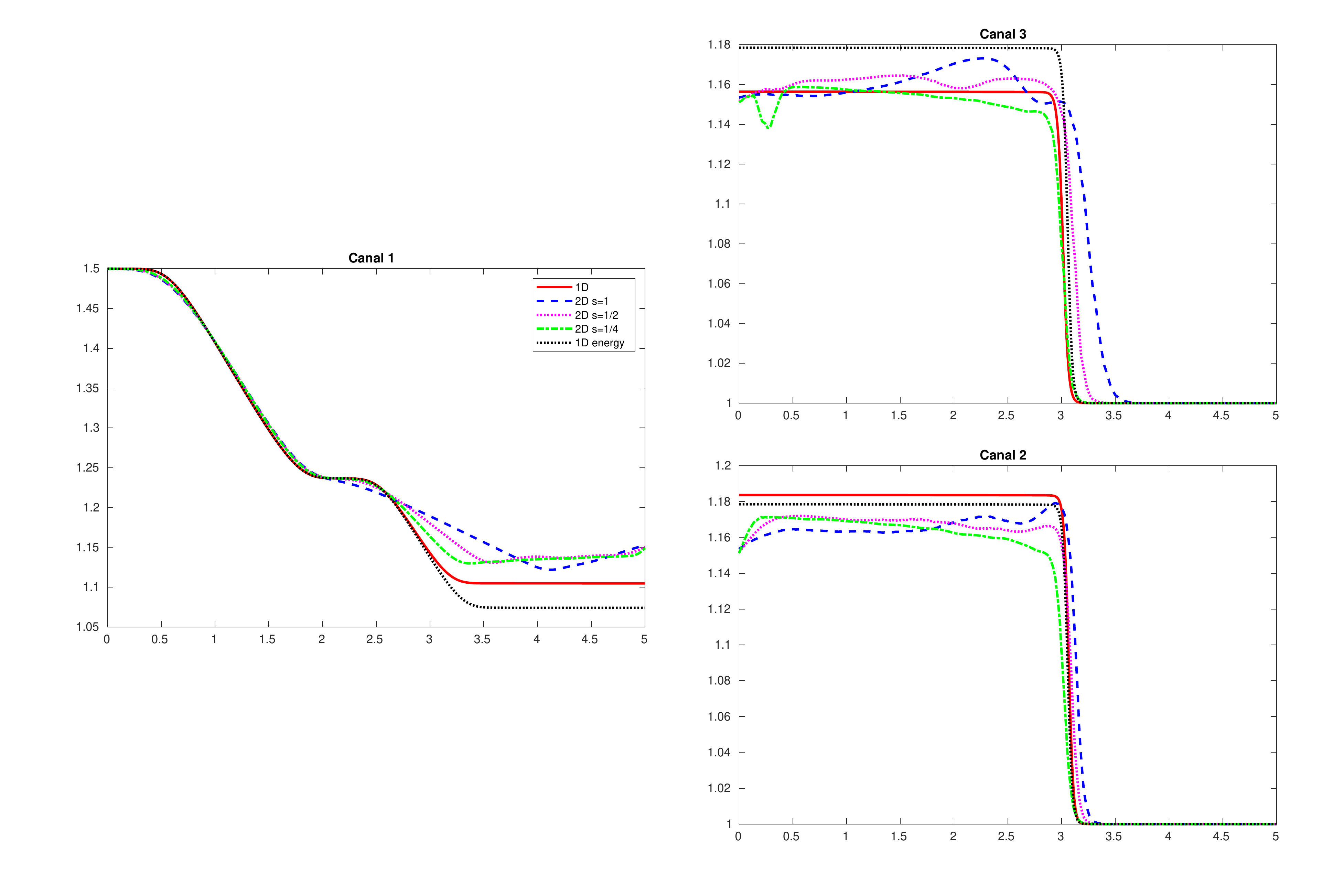}
  \end{overpic}
 }
 \\
  \subfloat[water velocity]
 {
 \begin{overpic}
 [width=1.0\textwidth]{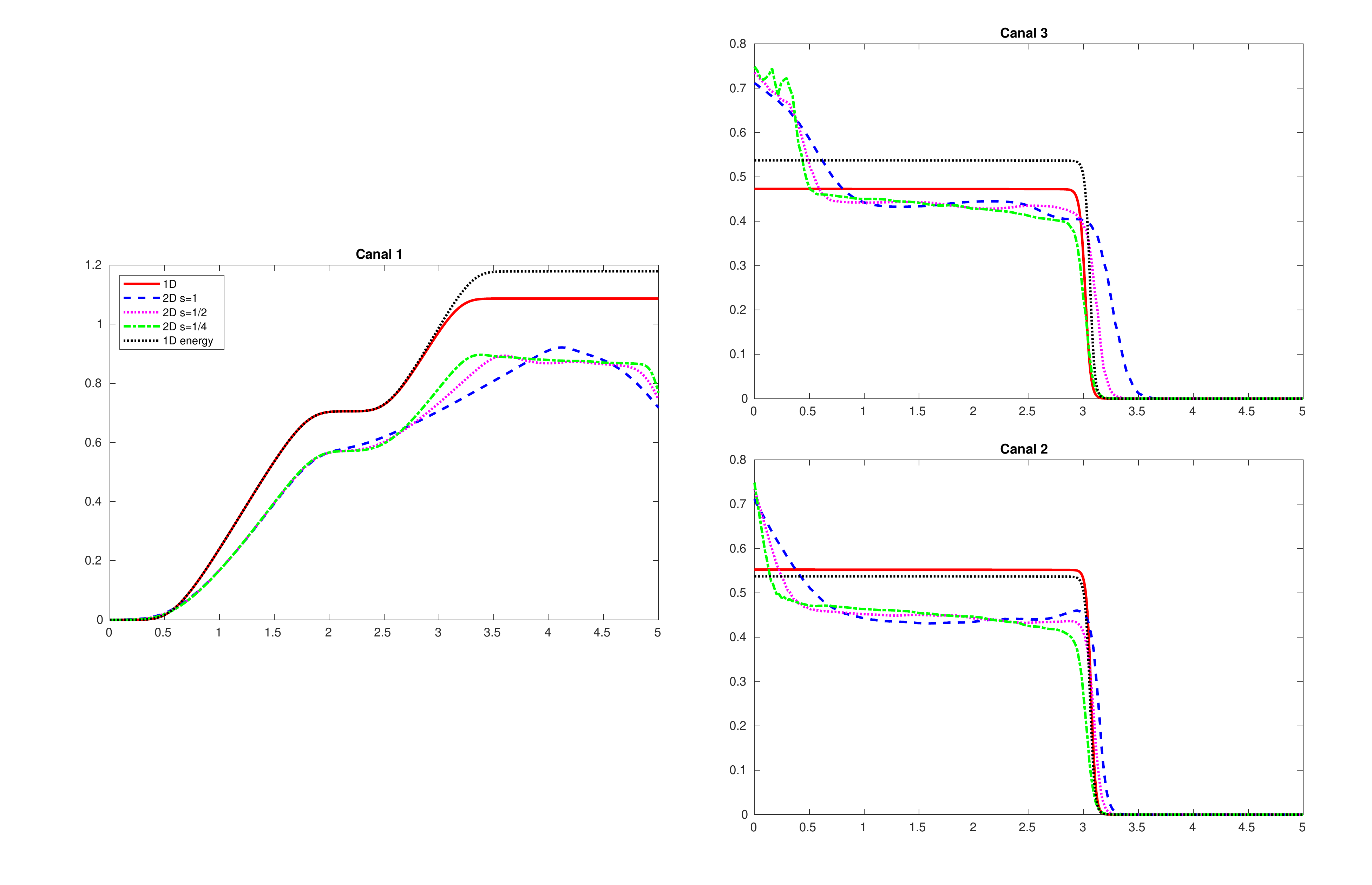}
  \end{overpic}
 }
\caption{
Test \ref{sec:test_1D2D}, \textit{diverging juntion}: comparison of 2D with 1D solutions at $T=1.5$: $\theta=\frac{\pi}{3}$ and $\phi=-\frac{\pi}{12}$, $L_1=L_2=L_3=5$, $s_1=s_2=s_3=s$ with $s=1,\ 0.5, \ 0.25$ (in the 2D case). The 1D solution with junction conditions \eqref{eq:Conservation}-\eqref{eq:cond_riemann_pb} is displayed in solid red, the 1D solution with equal energy condition in dotted black, the 2D solution with $s=1,\ 0.5,\ 0.25$ in dashed blue,  magenta and  green. 
Initial data are given at  Eq.\eqref{eq:initdata1} (1D case) and Eq.\eqref{eq:init2D_1} (2D case).
}\label{fig:confronto_red_sec}
\end{figure}

\paragraph{\it Diverging junction.} We fix $\theta=\pi/3$, $\phi=-\pi/12$ and initial states as in \eqref{eq:initdata1} for the 1D configuration and for the 2D system such that 
\begin{equation}\label{eq:init2D_1}
\begin{array}{l}
v_{x,y}(x,y,t=0)=0 \quad (x,y)\in\Omega,
\medskip\\
h(x,y,t=0) = \left\{\begin{array}{cl}
1.5 & (x,y)\in\Omega\cap\{ 0 \leq x\leq L_1/2\}
\\
1 & \mbox{otherwise}.
\end{array}\right.
\end{array}
\end{equation}
In Figure \ref{fig:confronto_red_sec}, the 1D solution with junction conditions \eqref{eq:Conservation}-\eqref{eq:cond_riemann_pb} is presented in solid red, the dotted black curve is the 1D solution assuming equal energy at the junction,  as in \cite{HS2013}.  We have three 2D solutions which are displayed in dashed blue,  magenta and  green for $s=1,\ 0.5,\ 0.25$ respectively. Decreasing $s$ we see that the 2D solutions converge towards the 1D wave front, both in the two shocks and the receding rarefaction. The only differences are observed in the flat intermediate states at the junction. The junction solver proposed in this work seems more accurate than the one with equal energy condition at the junction.
Some wiggles appearing in the 2D solution might be due to numerical artefacts at the interface between dry and wet states of the 2D code.

\medskip 

\paragraph{\it Merging junction.}
We fix $\theta=5\pi/12$, $\phi=-3\pi/8$ and initial data as in \eqref{eq:initdata2} for the 1D configuration while for the 2D system as \eqref{eq:init2D_1} with $h=1.5$ in the right half of channel 3 too, see Figure \ref{Fig:geom2D} on the right.

\begin{figure}[htbp!]
 \subfloat[water height]
 {
\begin{overpic}
 [width=0.9\textwidth]{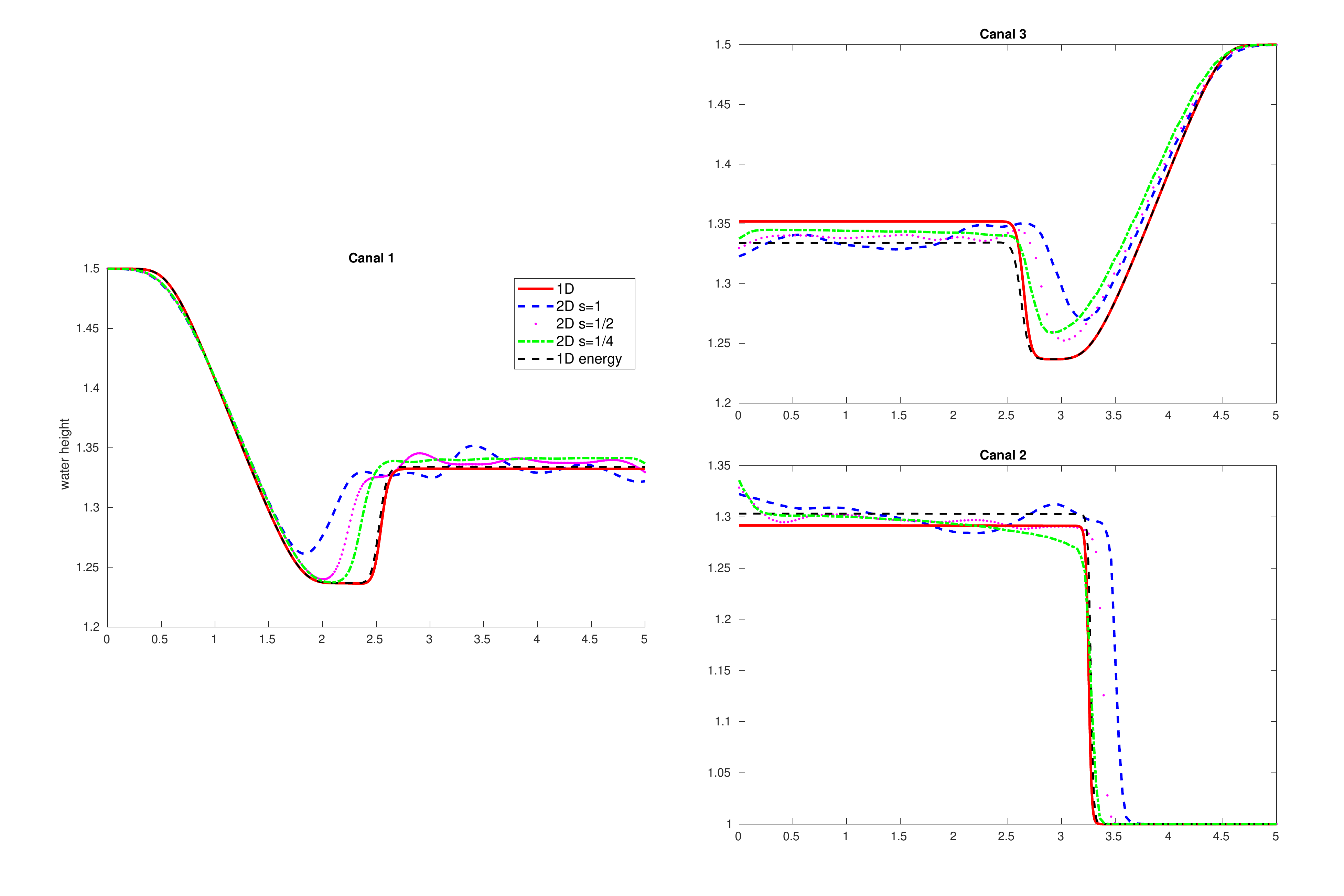}
  \end{overpic}
 }
 \\
  \subfloat[Module of water velocity]
 {
 \begin{overpic}
 [width=0.9\textwidth]{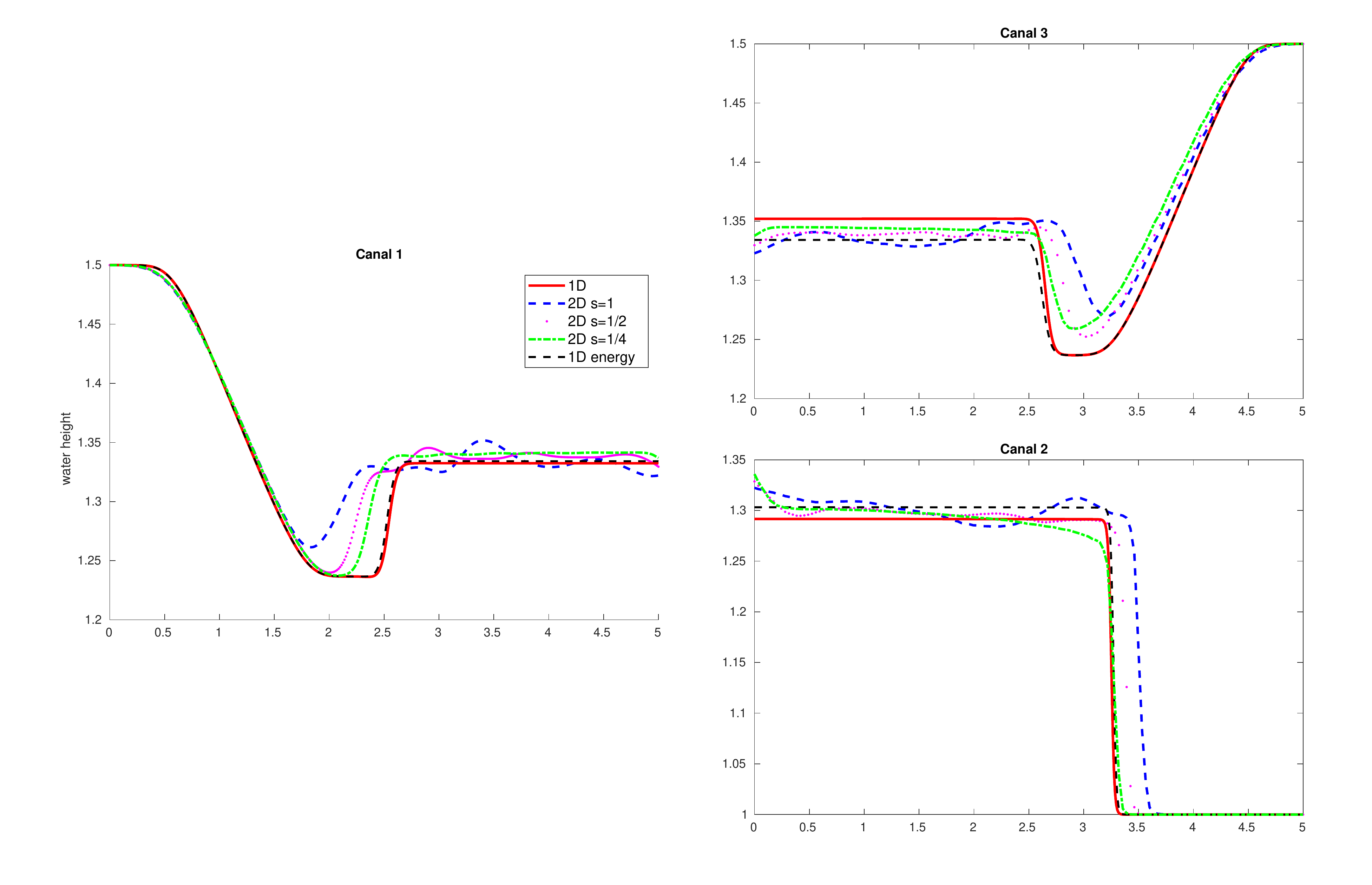}
  \end{overpic}
 }
\caption{Test \ref{sec:test_1D2D}, \textit{ merging junction}:  comparison of 2D with 1D solutions at $T=1.5$: $\theta=5\pi/12$, $\phi=-3\pi/8$, $L_1=L_2=L_3=5$, $s_1=s_2=s_3=s$ with $s=1,\ 0.5, \ 0.25$ (in the 2D case). The 1D solution with junction conditions \eqref{eq:Conservation}-\eqref{eq:cond_riemann_pb}  is displayed in red, the 1D solution with equal energy condition in dotted black, the 2D solution with $s=1,\ 0.5,\ 0.25$ in dashed blue, solid magenta and solid green. Initial data are given at  Eq.\eqref{eq:initdata2} (1D case) and Eq.\eqref{eq:init2D_2} (2D case).
}\label{fig:confronto_dueuno_2D}
\end{figure}

As before, in Figure \ref{fig:confronto_dueuno_2D}, we superpose the 1D and 2D curves. Again, our 1D solution is displayed in solid red and the dotted black curve represents the solution with equal energy condition at the junction, \cite{HS2013}. The 2D solutions are drawn in dashed  blue, magenta and green for $s=1,0.5,0.25$ respectively. Decreasing $s$,  we see that  the 2D solutions converge towards the 1D water front, as expected.

\section{Conclusions}
In this paper we have presented a numerical solver for one dimensional channels in a network. The solver is based on a finite volume scheme in each canal and coupling conditions at the junction are obtained with a single 2D element at the junction across which mass and the two components of momentum are conserved. This approach allows to take into account quite general geometries including the dependence on the angles with which the canals intersect at the junction and the sections of the canals. In this framework we can also include the construction of solvers for shallow-water problems along a channel with varying section.

The solver is based on the assumption that the flow across the junction is fluvial. Future work on this topic will be concentrated on the case of torrential flows and on the dependence of the bottom topography.

\subsection*{Acknowledgements}
This work was partly supported by MIUR (Ministry of University and Research) PRIN2017 project number 2017KKJP4X.


\bibliography{Junction}{}
\bibliographystyle{plain}

\end{document}